\documentclass[review]{elsarticle}
\usepackage[english,russian]{babel}
\usepackage{lineno,hyperref}
\modulolinenumbers[5]

\journal{Nonlinear  Analysis: Theory, Methods \& Application}
\usepackage[tbtags]{amsmath}
\usepackage{amsfonts,amssymb,mathrsfs,amscd,comment}
\newtheorem{Lm}{Lemma}[section]
\newtheorem{Ca}{Corollary}[section]
\newtheorem{Th}{Theorem}[section]

\newtheorem{Remark}{Remark}[section]

\newtheorem{Def}{Definition}[section]
\DeclareMathOperator*{\essinf}{ess\,inf}

\begin{document}

\begin{frontmatter}

\title{Traces of weighted Sobolev spaces with Muckenhoupt weight. The case $p=1$\tnoteref{mytitlenote}}
\tnotetext[mytitlenote]{
This work is supported by the RSF under a grant 14-50-00005}

\author{A.\,I.~Tyulenev}
\address{Steklov Mathematical Institute {\rm (}Russian Academy of Sciences{\rm )}}

\ead{tyulenev-math@yandex.ru, tyulenev@mi.ras.ru}

\numberwithin{equation}{section}

\begin{abstract}
A complete description of traces on  $\mathbb{R}^{n}$  of functions from the weighted Sobolev space
$W^{l}_{1}(\mathbb{R}^{n+1},\gamma)$, $l \in \mathbb{N}$, with weight  $\gamma \in A^{\rm loc}_{1}(\mathbb{R}^{n+1})$ is obtained.
In the case $l=1$ the proof of the trace theorem is based on a~special nonlinear algorithm for constructing a~system of tilings of the space~$\mathbb R^n$.
As the trace of the space $W^1_1(\mathbb R^{n+1},\gamma)$ we have the new function space $Z(\{\gamma_{k,m}\})$.

\end{abstract}


\end{frontmatter}
\linenumbers

\section{Introduction}
The problem of exact description of the trace space (on the boundary of a~domain) of a~weighted Sobolev space has a~long history.
A~short survey of the available results in this direction is given in~\cite{Tyu}. It is worth pointing out that,  since the appearance of
the pioneering work of Gagliardo~\cite{Gal} a~long time ago, it was only in \cite{Tyu},~\cite{Tyu2} that
a~complete description of the trace space on~$\mathbb{R}^{n}$ of the weighted Sobolev space $W^{l}_{p}(\mathbb{R}^{n+1},\gamma)$,  $p \in (1,\infty)$, with
weight $\gamma \in A_{p}^{\rm loc}(\mathbb{R}^{n+1})$ was obtained. The solution of this problem, with such a~high degree of generality, calls for the introduction
of new modifications of Besov-type spaces of variable smoothness and new machinery for studying thereof.
Thus, in the case $p \in (1,\infty)$ we have the result in the most general context possible at present.

The thing gets much worse in the case $p=1$.
Indeed, starting from 1957, as far as the author is aware, the bibliography on this subject lists only the papers \cite{Gin},~\cite{Tyu3}, which put forward
an exact description of the trace on~$\mathbb{R}^{n}$ by the weighted Sobolev space $W^{1}_{1}(\mathbb{R}^{n+1},\gamma)$.
However, a~weight~$\gamma$ in these papers was assumed to be a~model function depending only on some highlighted group of variables.
For example, in~\cite{Gin} it was assumed that a~weight~$\gamma$ depends only on the coordinate $x_{n+1}$ (note that in~\cite{Gin} a~more general multiweight case was considered,
when different derivatives in the Sobolev norm are integrated with different weights, but all weights depend on the same coordinate $x_{n+1}$), and in~\cite{Tyu3} it is assumed that
$\gamma \in A_{1}(\mathbb{R}^{n})$.

We note the recent paper \cite{Mi}, which contains many interesting results on the problem of exact description of the trace spaces of a weighted Sobolev spaces.
However, in this paper it was assumed that the weight $\gamma=\gamma(x_{1},...,x_{n+1})=|x_{n+1}|^{\alpha}$ with some constraints on the parameter~$\alpha$.

Of course, such a~lack of knowledge in the case $p=1$ is due to the great difficulty of the problem. Attempts to find the trace of the space
 $W^{l}_{1}(\mathbb{R}^{n+1},\gamma)$ with fairly general $\gamma$ involve considerable difficulties.

The machinery of \cite{Tyu}, \cite{Tyu2} may not in principle be applied in this setting,
because this approach  depends on the Muckenhoupt theorem on the boundedness  of the Hardy--\allowbreak Littlewood maximal operator in weighted Lebesgue spaces (this result fails for $p=1$,
see Ch.~5 of \cite{St} for the details).

However, in the `simple' nonlimiting case  $l > 1$ one eventually succeeds in
modifying the methods of \cite{Tyu},~\cite{Tyu2} (without having recourse to the Hardy--\allowbreak Littlewood maximal function!) to give a~complete description of
the trace space of the weighted Sobolev  space $W^{l}_{1}(\mathbb{R}^{n+1},\gamma)$ with weight $\gamma \in A_{1}^{\rm loc}(\mathbb{R}^{n+1})$ in terms of the
Besov spaces of variable smoothness that were introduced by the author~\cite{Tyu2}. We shall enlarge upon these results in~\S\,3.

The case $l=1$ presents the greatest challenge and calls for the development of a~refreshingly different method. The situation is aggravated by the fact that even in the case
$\gamma \equiv 1$ the extension operator from the trace space turns out to be nonlinear \cite{Gal},~\cite{Pe}. On the other hand in~\cite{Gin} it was shown that the extension operator
from the corresponding trace space is linear if $\lim_{x_{n+1} \to 0}\gamma(x_{n+1})=+\infty$ for a~continuous weight $\gamma=\gamma(x_{n+1})$.

\section{Basic notations and definitions}

As usual, $\mathbb{Z}_{+}$ and $\mathbb{N}$ will denote the set of all nonnegative and positive integers respectively. Also, $\mathbb{Z}^{n}$ is the linear space of vectors in~$\mathbb R^n$ with integer components.

Throughout we shall fix $n \in \mathbb{N}$, which will only be used to denote the dimension of the Euclidean space  $\mathbb{R}^{n}$.
A~point of the space $\mathbb{R}^{n}$ will be written as $x=(x_{1},\dotsc,x_{n})$, and a~point of the space  $\mathbb{R}^{n+1}$,
as the pair $(x,t)$ ($x \in \mathbb{R}^{n}$, $t \in \mathbb{R}$). The space $\mathbb{R}^{n}$ will be identified with the hyperplane $\mathbb{R}^{n} \times \{0\}$ of the space~$\mathbb{R}^{n+1}$.

The symbol $C$ will be used to denote (different) insignificant constants in various estimates. Sometimes, if it is required for purposes of
exposition, we shall indicate the parameters on which some or other constant depends.

The derivatives will be written in the multi-index notation: $
D^{\alpha}:=\frac{\partial^{|\alpha|}}{\partial x^{\alpha_{1}}_{1}\dots\partial x^{\alpha_{n}}_{n+1}}$, where
$\alpha$~is a~vector from~$\mathbb{Z}^{n+1}$ with nonnegative components ($\alpha=(\alpha_{1},\dotsc,\alpha_{n+1})$),
$|\alpha|:=|\alpha_{1}|+\dots+|\alpha_{n+1}|$.

Given a measurable subset $E$ of the space $\mathbb{R}^{d}$, $d =n, n+1$, we denote by $|E|$ the Lebesgue measure of~$E$, and by  $\chi_{E}$,
the characteristic function of~$E$.

By an open cube $Q$ in the space $\mathbb{R}^{d}$, $d = n,n+1$ (or simply a~cube, if the dimension of the ambient space is clear from the context)
we shall mean a~cube with sides parallel to coordinate axes. By $\overline{Q}$ we denote the closure of a~cube~$Q$ in the space~$\mathbb{R}^{d}$, $d=n,n+1$ which will be called a~closed cube.
By $r(Q)$ we denote the side length of~$Q$.

Given $k \in \mathbb{Z}_{+}$, $m=(m_{1},\dotsc,m_{d}) \in \mathbb{Z}^{d}$, we let
$Q_{k,m}:=\prod\limits_{i=1}^{d}(\frac{m_{i}}{2^{k}},\frac{m_{i}+1}{2^{k}})$
denote an open dyadic cube of rank~$k$ in the space $\mathbb{R}^{d}$, $d=n,n+1$.

Let $I:=\prod\limits_{i=1}^{n}(-1,1)$.

By a weight we shall imply a function $\gamma \in L^{\rm loc}_{1}(\mathbb{R}^{n+1})$ which is positive almost everywhere.

Next, for a measurable set $E \subset \mathbb{R}^{n+1}$ of positive measure and a~weight~$\gamma$, we define
$$
\gamma_{E}:=\frac{1}{|E|}\iint\limits_{E}\gamma(x,t)\,dx\,dt.
$$

It what follows we shall be concerned only with weights that locally satisfy the Muckenhoupt condition. Following \cite{Sa} we introduce

\begin{Def}\rm
\label{Def2.1}
We say that a weight $\gamma \in A^{\rm loc}_{1}(\mathbb{R}^{n+1})$ if
\begin{equation}
\label{eq2.1}
\gamma_{Q} \le C_{\gamma} \essinf_{(x,t) \in Q}\gamma(x,t)
\end{equation}
for any cube~$Q$ in~$\mathbb{R}^{n+1}$ of side length $r(Q) \le 1$. It is clear that $C_{\gamma} \geq 1$.
\end{Def}

\begin{Remark}\rm
\label{R2.1}
The next elementary observation will be used in the sequel. Let $Q$~be a~cube in $\mathbb{R}^{n+1}$ with side length $r(Q) \le 1$,
and let  $Q_{1} \subset Q$ be a~cube of halved side length. Using~\eqref{eq2.1}, we clearly have
\begin{gather} \notag
\essinf\limits_{(x,t) \in Q_{1}}\gamma(x,t) \le \gamma_{Q_{1}} \le \frac{|Q|}{|Q_{1}|}\gamma_{Q} \le C_{\gamma} 2^{n+1} \essinf\limits_{(x,t) \in Q}\gamma(x,t),\\
\iint\limits_{Q}\gamma(x,t)\,dx\,dt \le C_{\gamma}|Q| \essinf \limits_{(x,t) \in Q}\gamma(x,t) \le C_{\gamma}|Q| \essinf \limits_{(x,t) \in Q_{1}} \gamma(x,t) \le 2^{n+1}C_{\gamma}\iint\limits_{Q_{1}}\gamma(x,t)\,dx\,dt.
\label{eq2.2}
\end{gather}

From \eqref{eq2.2} one easily obtains that, for $k \in \mathbb{Z}_{+}$, $m,m' \in \mathbb{Z}^{n+1}$, $|m-m'|\le a$ ($a > 1$),
\begin{equation}
\label{eq2.3}
\begin{gathered}
\essinf \limits_{(x,t) \in Q_{k,m}}\gamma(x,t) \le C(C_{\gamma},n,a)  \essinf\limits_{(x,t) \in Q_{k,m'}}\gamma(x,t) ,\\
\iint\limits_{Q_{k,m}}\gamma(x,t)\,dx\,dt \le C(C_{\gamma},n,a)\iint\limits_{Q_{k,m'}}\gamma(x,t)\,dx\,dt.
\end{gathered}
\end{equation}

We give a detailed proof of only the first inequality in \eqref{eq2.3}, because the proof of the second one is similar.
Assume that for $k \in \mathbb{N}$ the cubes $Q_{k,m}$ and $Q_{k,m'}$ are disjoint, but have common boundary points.
Then there exists a~unique cube~$Q$ with twice greater side length containing both cubes $Q_{k,m}$ and $Q_{k,m'}$.
Besides, the side length of the cube~$Q$ is at most~$1$, because the side lengths of the cubes $Q_{k,m}$ and $Q_{k,m'}$ is at most $\frac{1}{2}$. From
\eqref{eq2.2} we get
\begin{equation}
\label{eq2.3'}
\begin{gathered}
\essinf \limits_{(x,t) \in Q_{k,m}}\gamma(x,t) \le 2^{n+1}C_{\gamma}\essinf \limits_{(x,t) \in Q}\gamma(x,t) \le 2^{n+1}C_{\gamma}\essinf \limits_{(x,t) \in Q_{k,m'}}\gamma(x,t).
\end{gathered}
\end{equation}

Now the first inequality in \eqref{eq2.3} with $k \geq 1$ clearly follows from estimate \eqref{eq2.3'}.

We also note that, for $k \in \mathbb{Z}_{+}$, $m \in \mathbb{Z}^{n}$,
\begin{equation}
\label{eq2.3''}
\begin{gathered}
\essinf \limits_{(x,t) \in Q_{k,m}}\gamma(x,t) = \min\limits_{Q_{k+1,m'} \subset Q_{k,m}}\quad \essinf \limits_{(x,t) \in Q_{k+1,m'}}\gamma(x,t).
\end{gathered}
\end{equation}

Now the required inequality (in the case $k=0$) easily follows from \eqref{eq2.3'}, \eqref{eq2.3''}.

Let $c > 1$, $\gamma \in A^{\rm loc}_{1}(\mathbb{R}^{n+1})$ and~$Q$ be an arbitrary cube in the space $\mathbb{R}^{n+1}$
with side length $r(Q) \le c$. Arguing as in the derivation of the first inequality of \eqref{eq2.3} and using \eqref{eq2.3'}, \eqref{eq2.3''}, this establishes
\begin{equation}
\label{eq2.4}
\begin{gathered}
\gamma_{Q} \le C_{1}(n,C_{\gamma},c) \sum\limits_{\substack {m \in \mathbb{Z}^{n+1}\\ Q_{0,m} \bigcap  Q \neq \emptyset}}\gamma_{Q_{0,m}} \le C_{2}(C_{\gamma},n,c)\sum\limits_{\substack {m \in \mathbb{Z}^{n+1}\\ Q_{0,m} \bigcap  Q \neq \emptyset}}\essinf \limits_{(x,t) \in Q_{0,m}}\gamma(x,t)\le \\
\le C_{3}(C_{\gamma},n,c) \min\limits_{\substack {m \in \mathbb{Z}^{n+1}\\ Q_{0,m} \bigcap  Q \neq \emptyset}} \quad \essinf \limits_{(x,t) \in Q_{0,m}}\gamma(x,t) \le C_{3}(C_{\gamma},n,c) \essinf \limits_{(x,t) \in Q}\gamma(x,t).
\end{gathered}
\end{equation}

Besides, for any cube~$Q$ in the space $\mathbb{R}^{n+1}$ with side length $r(Q) \le 1$,
\begin{equation}
\label{eq2.5}
\iint\limits_{cQ}\gamma(x,t)\,dx\,dt \le C(C_{\gamma},n,c)\iint\limits_{Q}\gamma(x,t)\,dx\,dt.
\end{equation}
\end{Remark}

The  proof of estimate \eqref{eq2.5} is similar to that of \eqref{eq2.4} and is based on the second inequality of~\eqref{eq2.3}.

\begin{Def}\rm
\label{Def2.2}
Assume that a weight $\gamma \in A^{\rm loc}_{1}(\mathbb{R}^{n+1})$, $l \in \mathbb{N}$, and $\Omega$~is a~domain in~$\mathbb{R}^{n+1}$.
By $W^{l}_{1}(\Omega,\gamma)$ we will denote the linear space of functions which are locally integrable  on $\Omega$  and have finite norm
\begin{equation}
\label{eq2.6}
\|f|W^{l}_{1}(\Omega,\gamma)\|:=\sum\limits_{|\alpha| \le l}\|\gamma D^{\alpha}f|L_{1}(\Omega)\|.
\end{equation}
\end{Def}

For $\gamma \equiv 1$ we shall write $W^{l}_{1}(\Omega)$ instead of $W^{l}_{1}(\Omega,1)$.

By $D^{\alpha}f$ in \eqref{eq2.6} we denote the (Sobolev) generalized derivatives of a~function~$f$ of order~$\alpha$ (see Ch.~1 of \cite{Bu} or Ch.~2 of~\cite{Zi}
for equivalent definitions and basic properties).

\begin{Remark}\rm
\label{R2.2}
Using \eqref{eq2.4} and H\"older's inequality, we see that  if the weight $\gamma \in A^{\rm loc}_{1}(\mathbb{R}^{n+1})$, $l \in \mathbb{N}$ and $\Omega$ is a~bounded domain
(in the space $\mathbb{R}^{n+1}$), then the space $W^{l}_{1}(\Omega,\gamma)$ is continuously embedded in the space~$W^{l}_{1}(\Omega)$
(with the embedding constant depending on $n$, $\operatorname{diam}\Omega$ and the constant $C_\gamma$).

\end{Remark}

The following fact will be frequently  useful. For completeness, we give the proof.

\begin{Lm}
\label{Lm2.1} Let $d \in \mathbb{N}$, $f \in L_{1}(\mathbb{R}^{d})$, and let $N \in \mathbb{N}$,
$\mathbb{R}^{d}=\bigcup\limits_{j=1}^{\infty}R_{j}$, where measurable sets~$R_{j}$ are such that any point $z \in \mathbb{R}^{d}$ lies in at most than~$N$ sets
from the family $\{R_{j}\}^{\infty}_{j=1}$. Then
$$
\sum\limits_{j=1}^{\infty}\int\limits_{R_{j}}|f(z)|\,dz \le N \|f|L_{1}(\mathbb{R}^{d})\|.
$$
\end{Lm}

\textbf{Proof.} From the hypotheses of the theorem  we see at once that
$$
\sum\limits_{j=1}^{\infty}\chi_{R_{j}}(z) \le N, \quad x \in \mathbb{R}^{d}.
$$
Hence, we have the estimate
\begin{gather*}
\sum\limits_{j=1}^{\infty}\int\limits_{R_{j}}|f(z)|\,dz=\sum\limits_{j=1}^{\infty}\int\limits_{\mathbb{R}^{n+1}}\chi_{R_{j}}(z)|f(z)|\,dz=
\\
=\int\limits_{\mathbb{R}^{n+1}}\sum\limits_{j=1}^{\infty}\chi_{R_{j}}(z)|f(z)|\,dz \le N \|f|L_{1}(\mathbb{R}^{n+1})\|.
\end{gather*}

\section{The nonlimiting case}

In this section we shall modify the methods of~\cite{Tyu},\cite{Tyu2} and give a~complete description of the trace space of the space
$W^{l}_{1}(\mathbb{R}^{n+1},\gamma)$ on the hyperplane $\mathbb{R}^{n}$ under the condition that $\gamma \in A^{\rm loc}_{1}(\mathbb{R}^{n+1})$, $l > 1$.
Until the end of this section,  $Q_{k,m}$ ($k \in \mathbb{Z}_{+}$, $m \in \mathbb{Z}^{n}$) will denote dyadic cubes of rank~$k$ in the space $\mathbb{R}^{n}$.

For the rest of this section we fix a parameter $l \in \mathbb{N}$, $l > 1$, and a~weight $\gamma \in A^{\rm loc}_{1}(\mathbb{R}^{n+1})$. Next, we set
$$
\gamma_{k,m}:=2^{kl}\iint\limits_{Q_{k,m} \times (0,2^{-k})}\gamma(x,t)\,dx\,dt, \qquad (k,m) \in \mathbb{Z}_{+} \times \mathbb{Z}^{n}.
$$

\begin{Remark}\rm
\label{R3.1}
From \eqref{eq2.3} it clearly follows that
$\gamma_{k,m} \le C_{1}(C_{\gamma},n,c) \gamma_{k,m'}$  for  $k \in \mathbb{Z}_{+}$, $m,m' \in \mathbb{Z}^{n}$, $|m-m'| \le c$ ($c \geq 1$). Furthermore from \eqref{eq2.2} we have
$\gamma_{k+1,m'} \le C_{2}(C_{\gamma},n,l)\gamma_{k,m}$, $\gamma_{k,m} \le C_{3}(C_{\gamma},n,l)\gamma_{k+1,m'}$
for  $k \in \mathbb{Z}_{+}$, $m,m' \in \mathbb{Z}^{n}$  and  $Q_{k+1,m'} \subset Q_{k,m}$.

\end{Remark}

For further purposes we shall need the definition of the Besov-type space of variable smoothness. Actually, we give a~particular case of Definition 2.5 of~\cite{Tyu},
because we shall not need the whole range of the parameters (and such general assumptions on the variable smoothness).

Given a measurable function~$\varphi$  and $x,h \in \mathbb{R}^{n}$, we define
$\Delta^{l}(h)\varphi(x):=\sum\limits_{i=0}^{l}(-1)^{i}C^{i}_{l}\varphi(x+ih)$. Next, for a~function $\varphi \in L^{\rm loc}_{1}(\mathbb{R}^{n})$ and
a~cube~$Q$ in the space~$\mathbb{R}^{n}$, we set
$$
\delta^{l}(Q)\varphi:=\frac{1}{|Q|^{2}}\int\limits_{r(Q)I}\int\limits_{Q}|\Delta^{l}(h)\varphi(x)|\,dxdh.
$$

By  $E^{l}(Q)\varphi$ we shall denote the best $L_{1}(Q)$-approximation to a~function $\varphi \in L^{\rm loc}_{1}(\mathbb{R}^{n})$ on a~cube $Q$ by polynomials
of degree~$<l$.

\begin{Def} \rm
\label{Def3.1} By $\widetilde{B}^{l}(\mathbb{R}^{n},\{\gamma_{k,m}\})$ we shall denote the Besov-type space of variable smoothness $\{\gamma_{k,m}\}$ equipped with the norm
\begin{equation}
\label{eq3.1}
\left\|\varphi|\widetilde{B}^{l}(\mathbb{R}^{n},\{\gamma_{k,m}\})\right\|:=\sum\limits_{k=1}^{\infty}\sum\limits_{m \in \mathbb{Z}^{n}}
\gamma_{k,m}\delta^{l}(Q_{k,m})\varphi+\sum\limits_{m \in \mathbb{Z}^{n}}\gamma_{0,m}\|\varphi|L_{1}( Q_{0,m})\|.
\end{equation}
\end{Def}

\begin{Remark}\rm
\label{R3.2}
According to \cite{Tyu},  for $c > 1$,
$$
\left\|\varphi|\widetilde{B}^{l}(\mathbb{R}^{n},\{\gamma_{k,m}\})\right\| \sim \sum\limits_{k=1}^{\infty}\sum\limits_{m \in \mathbb{Z}^{n}}
\gamma_{k,m}2^{kn}E^{l}(cQ_{k,m})\varphi+\sum\limits_{m \in \mathbb{Z}^{n}}\gamma_{0,m}\|\varphi|L_{1}( Q_{0,m})\|.
$$
\end{Remark}

\begin{Def}\rm
\label{Def3.2}
A function $\varphi \in L^{\rm loc}_{1}(\mathbb{R}^{n})$ is said to be the trace of a~function $f \in L^{\rm loc}_{1}(\mathbb{R}^{n+1})$ on the hyperplane $\mathbb{R}^{n}$ (written $\operatorname{tr}\left|_{t=0}f=\varphi\right.$) if, for any open~$Q$ (in the space $\mathbb{R}^{n}$),
\begin{equation}
\label{eq3.2}
\int\limits_{Q}|\varphi(x)-f(x,t)|\,dx  \to 0, \qquad  t \to 0.
\end{equation}
\end{Def}

Let $E \subset L^{\rm loc}_{1}(\mathbb{R}^{n+1})$ be the linear space of functions~$f$ that have the trace on the hyperplane~$\mathbb R^n$.
In what follows, by $\operatorname{Tr}$ we shall denote the linear operator $\operatorname{Tr}:E \to L^{\rm loc}_{1}(\mathbb{R}^{n})$
defined by $\operatorname{Tr}[f]=\operatorname{tr}\left|_{t=0}f=\varphi\right.$.

\begin{Th}
\label{Th3.1}
The linear operator $\operatorname{Tr}:W^{l}_{1}(\mathbb{R}^{n+1},\gamma) \to \widetilde{B}^{l}(\mathbb{R}^{n},\{\gamma_{k,m}\})$ is bounded.
Moreover, there exists a~bounded linear operator $\operatorname{Ext}:\widetilde{B}^{l}(\mathbb{R}^{n},\{\gamma_{k,m}\}) \to W^{l}_{1}(\mathbb{R}^{n+1},\gamma)$
such that $\operatorname{Tr} \circ \operatorname{Ext} = \operatorname{Id}$ on the space  $\widetilde{B}^{l}(\mathbb{R}^{n},\{\gamma_{k,m}\})$.
\end{Th}

\textbf{Proof.} \textit{Step $1$}. Assume that a function $f \in W^{l}_{1}(\mathbb{R}^{n+1},\gamma)$. Then Remark~\ref{R2.2} and Theorem~2 of \S\,5.2 of~\cite{Bu} show that
the function~$f$ has the trace (which we denote by~$\varphi$) on~$\mathbb{R}^{n}$ and moreover

\begin{equation}
\label{eq3.3}
\varphi(x)=\lim\limits_{t \to +0}f(x,t) \qquad \hbox{for almost all } x \in \mathbb{R}^{n}.
\end{equation}

Let us prove the estimate
\begin{equation}
\label{eq3.4}
\|\varphi|\widetilde{B}^{l}(\mathbb{R}^{n},\{\gamma_{k,m}\})\| \le C \|f|W^{l}_{1}(\mathbb{R}^{n+1},\gamma)\|,
\end{equation}
where the constant $C > 0$ is independent of the function~$f$.

Let  $m \in \mathbb{Z}^{n}$ be fixed. By Remark~\ref{R2.2}, using the definition of the (Sobolev) generalized derivative of~$f$, and \eqref{eq3.3} we have

\begin{equation}
\label{eq3.5}
f(x,t)-\varphi(x)=\int\limits_{0}^{t}D_{t}f(x,\tau)\,d\tau \qquad \hbox{for almost all } x \in \mathbb{R}^{n}.
\end{equation}

Using \eqref{eq3.5}, \eqref{eq2.1} this gives
\begin{equation}
\label{eq3.6}
\begin{gathered}
\gamma_{0,m}\int\limits_{Q_{0,m}}|\varphi(x)|\,dx \le
\gamma_{0,m}\int\limits_{Q_{0,m}}\int\limits_{0}^{1}|f(x,\tau)-\varphi(x)|+|f(x,\tau)|\,dx\,d\tau \le \\ \le
C_{\gamma} \left(\int\limits_{Q_{0,m}}\int\limits_{0}^{1}\gamma(x,\tau)|D_{t}f(x,\tau)|\,dx\,d\tau + \int\limits_{Q_{0,m}}\int\limits_{0}^{1}\gamma(x,\tau)|f(x,\tau)|\,d\tau dx \right).
\end{gathered}
\end{equation}

Summing estimate  \eqref{eq3.6} over all  $m \in \mathbb{Z}^{n}$, we see that
\begin{equation}
\label{eq3.7}
\sum\limits_{m \in \mathbb{Z}^{n}}\gamma_{0,m}\int\limits_{Q_{0,m}}|\varphi(x)|\,dx \le C_{\gamma} \|f|W^{l}_{1}(\mathbb{R}^{n+1},\gamma)\|.
\end{equation}

From estimates (3.4) of Lemma 3.1 of \cite{Tyu2}, we get the estimate
\begin{equation}
\label{eq3.8}
\delta^{l}(Q_{k,m})\varphi \le  C 2^{kn} \int\limits_{C_{1}Q_{k,m}}\int\limits^{C_{2}2^{-k}}_{0}\sum\limits_{|\alpha|=l}t^{l-1}|D^{\alpha}f(x,t)|\,dt\,dx,
\end{equation}
in which the constants $C,C_{1},C_{2}$ depend only on $l,n$.

An application of \eqref{eq2.1} gives
\begin{equation}
\label{eq3.9}
\begin{gathered}
\sum\limits_{j=1}^{k}\sum\limits_{\substack {m' \in \mathbb{Z}^{n}\\ Q_{k,m} \subset Q_{j,m'}}}2^{jn}\gamma_{j,m'} = \sum\limits_{j=1}^{k}\sum\limits_{\substack {m' \in \mathbb{Z}^{n}\\ Q_{k,m} \subset Q_{j,m'}}}2^{j(l-1)}2^{j(n+1)}\iint\limits_{Q_{j,m'}\times(0,2^{-j})}\gamma(x,t)\,dx\,dt \le\\
\le C_{\gamma} \sum\limits_{j=1}^{k}\sum\limits_{\substack {m' \in \mathbb{Z}^{n}\\ Q_{k,m} \subset Q_{j,m'}}}2^{j(l-1)}\essinf_{(x,t) \times Q_{j,m'}\times(0,2^{-j})}\gamma(x,t) \le C_{\gamma} 2^{(k+1)(l-1)}\essinf_{(x,t) \in  Q_{k,m}\times(0,2^{-k})}\gamma(x,t).
\end{gathered}
\end{equation}

For brevity, we put $g(x,t):=\sum\limits_{|\alpha|=l}t^{l-1}|D^{\alpha}f(x,t)|$ for $(x,t) \in \mathbb{R}^{n+1}$.

We fix the smallest number  $k_{0} \in \mathbb{Z}_{+}$ for which $2^{k_{0}} > C_{2}$ (the constant $C_{2}$ is the same as on the right of~\eqref{eq3.8}).

The following equality holds for $k \in \mathbb{Z}_{+}$, $m \in \mathbb{Z}^{n}$
\begin{gather}
\int\limits_{Q_{k,m}}\int\limits_{0}^{2^{-(k-k_{0})}}g(x,t)\,dt dx = \sum\limits_{j=k}^{\infty}\sum\limits_{\substack {m' \in \mathbb{Z}^{n}\\ Q_{j,m'} \subset Q_{k,m}}}\int\limits_{Q_{j,m'}}\int\limits_{2^{-j-1}}^{2^{-j}}g(x,t)\,dt dx + \sum\limits_{j=k-k_{0}}^{k-1}\int\limits_{Q_{k,m}}\int\limits_{2^{-j-1}}^{2^{-j}}g(x,t)\,dt dx
\label{eq3.10}
\end{gather}

Clearly, for $k \in \mathbb{N}$, $m \in \mathbb{Z}^{n}$,
$$
\delta^{l}(Q_{k,m})\varphi \le  2^{2kn}\|\varphi|L_{1}(C(n,l)Q_{k,m})\|.
$$

Hence, using the finite overlapping multiplicity of the sets $C(n,l)Q_{k,m}$ with $k \in \{1,..,k_{0}\}$, $m \in \mathbb{Z}^{n}$, and using  Remark~\ref{R3.1}, it
follows from \eqref{eq3.7} that
\begin{gather}
S_{1}:=\sum\limits_{k=1}^{k_{0}}\sum\limits_{m \in \mathbb{Z}^{n}}\gamma_{k,m} \delta^{l}(Q_{k,m})\varphi \le C(n,l,k_{0},C_{\gamma})\sum\limits_{m \in \mathbb{Z}^{n}}\gamma_{0,m}\|\varphi|L_{1}(Q_{0,m})\| \le \notag \\
\le C(n,l,k_{0},C_{\gamma}) \|f|W^{l}_{1}(\mathbb{R}^{n+1},\gamma)\|.
\label{eq3.11}
\end{gather}

The sets $C_{1}Q_{k,m} \times (0,2^{-(k-k_{0})})$
have finite (depending only on $n$,~$l$) overlapping multiplicity (when index $k \in \mathbb{Z}_{+}$ is fixed and $m \in \mathbb{Z}^{n}$ variable), and hence, using  Remark~\ref{R3.1} it follows from \eqref{eq3.8}  that

\begin{gather}
\allowdisplaybreaks
S_{2}:=\sum\limits_{k=k_{0}+1}^{\infty}\sum\limits_{m \in \mathbb{Z}^{n}}
\gamma_{k,m}\delta^{l}(Q_{k,m})\varphi \le  C \sum\limits_{k=k_{0}+1}^{\infty}\sum\limits_{m \in \mathbb{Z}^{n}}
2^{kn} \gamma_{k,m}\int\limits_{Q_{k,m}}\int\limits^{2^{-(k-k_{0})}}_{0}g(x,t)\,dx\,dt
\label{eq3.12}
\end{gather}

Substituting \eqref{eq3.10} into the right-hand side of \eqref{eq3.12}, changing the order of summation  (in $k$ and~$j$), using Remark~\ref{R3.1} and estimate \eqref{eq3.9}, we obtain
\begin{gather}
S_{2} \le  C \sum\limits_{j=1}^{\infty}\sum\limits_{m' \in \mathbb{Z}^{n}}\int\limits_{Q_{j,m'}}\int\limits_{2^{-j-1}}^{2^{-j}}g(x,t)\,dt dx \left(\sum\limits_{k=1}^{j}\sum\limits_{\substack {m \in \mathbb{Z}^{n}\\ Q_{j,m'} \subset Q_{k,m}}}2^{kn}\gamma_{k,m} + \sum\limits_{k=j+1}^{j+k_{0}}\sum\limits_{\substack {m \in \mathbb{Z}^{n}\\ Q_{j,m'} \supset Q_{k,m}}}2^{kn}\gamma_{k,m}\right) \le \notag \\
\le C \sum\limits_{j=1}^{\infty}\sum\limits_{m \in \mathbb{Z}^{n}}2^{j(l-1)}\inf\limits_{(x,t) \times Q_{k,m}\times(0,2^{-j})}\gamma(x,t)\iint\limits_{Q_{k,m}\times(2^{-j-1},2^{-j})}g(x,t)\,dx\,dt \le  \notag \\
\le C\sum\limits_{j=1}^{\infty}\sum\limits_{m \in \mathbb{Z}^{n}}\iint\limits_{Q_{j,m}\times(2^{-j-1},2^{-j})}\gamma(x,t)\sum\limits_{|\alpha|=l}|D^{\alpha}f(x,t)|\,dx\,dt \le C \|f|W^{l}_{1}(\mathbb{R}^{n+1},\gamma)\|,
\label{eq3.13}
\end{gather}

in which the constant $C > 0$ depends only on $n,l,C_{\gamma},k_{0}$.

Now estimate \eqref{eq3.3} follows from \eqref{eq3.7}, \eqref{eq3.11}, \eqref{eq3.13}.

\textit{Step $2$}.
The construction of the extension operator  $\operatorname{Ext}:\widetilde{B}^{l}(\mathbb{R}^{n},\{\gamma_{k,m}\}) \to W^{l}_{1}(\mathbb{R}^{n+1},\gamma)$ and
the proof of its boundedness require minor modifications of the Step~2 in the proof of Theorem~3.1 from~\cite{Tyu}.
However, for the completeness of exposition we give the detailed proof.
So let
$\{\psi_k\}_{k=0}^\infty$  be a~partition of unity for the ball~$B^1:=(-1,1)$; that is, $\psi_k(t)\ge 0$ for $k\in \mathbb Z_+$, $t \in B^1$ and
$\sum_{k=0}^\infty\psi_k(t)=1$ for $t\in B^1$. Besides,
\begin{gather*}
\psi_0\in C^\infty\biggl(B^1 \setminus\frac{1}{2}B^1\biggr),
\qquad
\psi_k\in C_0^\infty\biggl(\frac{1}{2^{k-1}}B^1
\setminus\frac{1}{2^{k+1}}B^1\biggr)
\quad\text{for }\
k\in\mathbb N,
\\
|D^\beta\psi_k(t)|
\le C_{1} 2^{k|\beta|}
\quad\text{for}\quad
t\in B^1,
\quad
k\in\mathbb {Z}_{+}.
\end{gather*}

Assume that, for any $k\in\mathbb Z_+$,
only two functions $\psi_k$ and~$\psi_{k+1}$ 
do not vanish on the set $2^{-k}B^1 \setminus 2^{-k-1}B^1$. Hence,
$D^\beta\psi_k(t)=-D^\beta\psi_{k+1}(t)$
for $t\in 2^{-k}B^1 \setminus 2^{-k-1}B^1$.
The existence of a~sequence
$\{\psi_k\}_{k=0}^\infty$
 with the above properties may be proved as it was done, for example,
in \S\,4.5 of the book~\cite{Bu} in the proof of the trace theorem for unweighted Sobolev spaces.

We set
$$
f(x,t):=\sum_{k=1}^\infty\psi_k(t)
E_{2^{-k}}[\varphi](x)\quad
\text{for }\
(x,t)\in\mathbb R^{n+1},
$$
where, the operator $E_\varepsilon$  
(with $\varepsilon>0$) is defined as follows.

For a~function  $\varphi\in L^{\mathrm{loc}}_1(\mathbb R^n)$
we set
\begin{equation}
E_\varepsilon[\varphi](x)
:=\frac{1}{\varepsilon^{2n}}
\sum_{j=1}^l\mu_j\int_{\mathbb R^n}
\Theta\biggl(\frac{y-x}{\varepsilon}\biggr)
\int_{\mathbb R^n}\Theta
\biggl(\frac{z-y}{j\varepsilon}\biggr)
\varphi(z)\,dz\,dy, \quad  x\in\mathbb R^n.
\label{eq3.14}
\end{equation}

We shall not write down precise expressions for the constants  $\mu_{j}$ and the function~$\Theta$ from~\eqref{eq3.14},
which may be found in the authors' paper ~\cite{Tyu2} (Section~4).
The only important thing for us is that $\Theta \in C^{\infty}_{0}(B^{n})$, $\int\limits_{\mathbb{R}^{n}}\Theta(x)\,dx = 1$.

The following estimates are also of great value for us. Their proofs may be found in the authors' papers \cite{Tyu} (Lemma 3.1) and \cite{Tyu2} (Lemma~4.2).

A multi-index $\alpha=(\alpha_1,\dots,\alpha_{n+1})$ will be written as $(\alpha',\alpha_{n+1})$.

For any number $\varepsilon>0$, a~multi-\allowbreak index~$\alpha'=(\alpha_{1},..,\alpha_{n})$, $|\alpha'|=l$, and
$x\in\mathbb R^n$,
\begin{equation}
|D^{\alpha'}E_\varepsilon[\varphi](x)|
\le\frac{C}{\varepsilon^l}\mspace{2mu}
\delta^l(x+\varepsilon I)\varphi.
\label{eq3.15}
\end{equation}

Moreover, for any numbers $0<\varepsilon_1<\varepsilon_2$
a multi-index~$\beta'=(\beta_{1},..,\beta_{n})$, and $x\in\mathbb R^n$,
\begin{equation}
|D^{\beta'} E_{\varepsilon_1}[\varphi](x)
-D^{\beta'} E_{\varepsilon_2}[\varphi](x)|
\le C\int_{\varepsilon_1}^{\varepsilon_2}
\frac{1}{t^{1+|\beta|}}\mspace{2mu}
\delta^l(x+t I)\varphi\,dt.
\label{eq3.16}
\end{equation}

For later purposes we note that by the construction the function~$f$ vanishes on the set $\mathbb{R}^{n} \times (\mathbb{R} \setminus B^1)$.

We set $\Xi_{k,m}^{1,n}:=Q_{k,m} \times (2^{-k}B^{1}\setminus 2^{-k-1}B^{1})$ for $k \in \mathbb{Z}_{+}$, $m \in \mathbb{Z}^{n}$.

Clearly,
\begin{align*}
&\iint_{\mathbb R^n\times\frac{1}{2}B^1}
\gamma(x,t)
\biggl\{\sum_{|\alpha|=l,\,\alpha_{n+1}=0}
|D^{\alpha}f(x,t)|
+\sum_{|\alpha|=l,\,\alpha_{n+1}>0}
|D^\alpha f(x,t)|\biggr\}\,dx\,dt
\\
&\qquad=\sum_{k=1}^\infty
\sum_{m\in\mathbb Z^n}
\iint_{\Xi_{k,m}^{1,n}}\gamma(x,t)
\\
&\qquad\qquad\times
\biggl\{\sum_{|\alpha|=l,\,\alpha_{n+1}=0}
|D^\alpha f(x,t)|
+\sum_{|\alpha|=l,\,\alpha_{n+1}>0}|
D^\alpha f(x,t)|\biggr\}\,dx\,dt.
\end{align*}

Taking into account properties of the functions $\psi_{k}$ and applying estimate~\eqref{eq3.16}, we see that
{\allowdisplaybreaks
\begin{align}
\notag
&\sum_{|\alpha|=l,\,\alpha_{n+1}>0}
\iint_{\Xi_{k,m}^{1,n}}\gamma(x,t)
|D^{\alpha}f(x,t)|\,dx\,dt
=\sum_{|\alpha|=l,\,\alpha_{n+1}>0}
\iint_{\Xi_{k,m}^{1,n}}\gamma(x,t)
\\
&\qquad\qquad\times\bigl|D^{\alpha_{n+1}}\psi_k(t)
D^{\alpha'}E_{2^{-k}}\varphi(x)
+D^{\alpha_{n+1}}\psi_{k+1}(t)
D^{\alpha'}E_{2^{-(k+1)}}\varphi(x)\bigr|\,dx\,dt
\nonumber
\\
&\qquad\le\sum_{|\alpha'|=l-\alpha_{n+1}}2^{k\alpha_{n+1}}
\iint_{\Xi_{k,m}^{1,n}}\gamma(x,t)
|D^{\alpha'}E_{2^{-k}}\varphi(x)
-D^{\alpha'}E_{2^{-(k+1)}}\varphi(x)|\,dx\,dt
\nonumber
\\
&\qquad\le C2^{2nk}\gamma_{k,m}
\biggl[\int_{\widetilde CQ_{k,m}}
\int_{{I}/{2^k}}
|\Delta^l(h)\varphi(z)|\,dh\,dz\biggr]\,dx\,dt
\notag
\\
&\qquad\qquad\qquad\qquad\qquad\qquad
\text{for}\quad
k\in\mathbb N,\quad
m\in\mathbb Z^n.
\label{eq3.17}
\end{align}
}

The constant $\widetilde C \geq 1$, which is the dilation coefficients of the cubes  $Q_{k,m}$, depends only on $l$, $n$ and the diameter of the support of the function~$\Theta$ from~\eqref{eq3.14}.

Similarly, it follows from \eqref{eq3.15} that
\begin{align}
&\sum_{|\alpha'|=l}
\iint_{\Xi_{k,m}^{1,n}}
\gamma(x,t)|D^{\alpha'}f(x,t)|\,dx\,dt
\nonumber
\\
&\qquad\le C\sum_{|\alpha'|=l}
\iint_{\Xi_{k,m}^{1,n}}\gamma(x,t)\max
\bigl\{|D^{\alpha'}E_{2^{-k}}\varphi(x)|,
|D^{\alpha'}E_{2^{-k-1}}\varphi(x)|\bigr\}\,dx\,dt
\nonumber
\\
&\qquad\le C2^{2nk}\gamma_{k,m}
\biggl[\int_{\widetilde{C} Q_{k,m}}
\int_{{I}/{2^k}}
|\Delta^l(h)\varphi(z)|\,dh\,dz\biggr]\quad
\text{for}\quad
k\in\mathbb {N},\quad
m\in\mathbb {Z}^{n}.
\label{eq3.18}
\end{align}

Using the definition of the function $f$, we have, for  $|\alpha|=l$,
\begin{align}
&\sum_{|\alpha|=l}
\iint_{\mathbb{R}^{n+1} \setminus (\mathbb{R}^{n} \times \frac{1}{2}B^{1})}
\gamma(x,t)|D^\alpha f(x,t)|\,dx\,dt \le
\nonumber
\\
&\qquad\le C\sum_{m\in\mathbb Z^n}\gamma_{0,m}
\|\varphi\mid L_1(\widetilde CQ_{0,m})\|
\le C\sum_{m\in\mathbb Z^n}\gamma_{0,m}
\|\varphi\mid L_1(Q^n_{0,m})\|,
\label{eq3.19}
\end{align}
since the cubes $\widetilde{C} Q^{n}_{k,m}$ have finite overlapping multiplicity  (the constant  $\widetilde{C}$ is the same as in \eqref{eq3.17}).

Hence, summing up estimates \eqref{eq3.17}, \eqref{eq3.18} in $k$ and~$m$, taking into account that the
cubes $ \widetilde CQ_{k,m}^n$ have finite overlapping multiplicity (with fixed $k\in\mathbb N$ and variable $m\in\mathbb Z^n$), and employing estimate \eqref{eq3.19}, this gives
\begin{equation}
\sum_{|\alpha|=l}\|D^\alpha f
\mid L_1(\mathbb R^{n+1},\gamma)\|
\le C\|\varphi\mid\widetilde{B}^l
(\mathbb R^n,\{\gamma_{k,m}\})\|.
\label{eq3.20}
\end{equation}

To estimate the generalized derivatives $D^{\alpha}f$ for $|\alpha| < l$ we write, for each $(x,t) \in \mathbb{R}^{n} \times B^{1}$, the integral representation
of the function $D^{\alpha}f$ in the cone (see \S\,3.4, \cite{Bu}),
$$
V(x,t)=\biggl\{(x,t)(1-\xi)+\xi(x',t')\mid \xi \in[0,1],\,
(x',t')\in\frac{1}{2}B^{n+1}(x,t+3)\biggr\}
$$
(here $\frac{1}{2}B^{n+1}(x,t+3)$ is the ball of radius
$\frac{1}{2}$ centred at $(x,t+3)$), and use Remark~16 of \S\,3.5 in~\cite{Bu}.

Let $|\alpha|<l$. Since $f(x,t)=0$ for $|t| > 1$, we have
$$
|D^\alpha f(x,t)|\le C\sum_{|\beta|=l}
\iint_{(x,0)+(I \times B^1)}
|D^\beta f(\widetilde x,\widetilde t)|\,
d\widetilde x\,d\widetilde t\quad
\text{for }\
(x,t)\in\mathbb R^n\times B^1.
$$

Hence employing H\"older's inequality and \eqref{eq2.1}, \eqref{eq2.3}, we obtain, for $m \in \mathbb{Z}^{n}$, $|\alpha| < l$,
{\allowdisplaybreaks
\begin{align}
&\iint_{Q_{0,m}\times B^1}\gamma(x,t)
|D^\alpha f(x,t)|\,dx\,dt
\nonumber
\\
&\qquad\le C\sum_{|\beta|=l}
\bigl[\iint\limits_{\widetilde CQ_{0,m}\times B^1}\gamma(x,t)\,dxdt\bigr]
\bigl[\essinf_{\widetilde CQ_{0,m}\times B^1}\gamma(x,t)]^{-1}
\nonumber
\\
&\qquad\qquad\times
\iint_{\widetilde CQ_{0,m}\times B^1}\gamma(x,t)
|D^\beta f(x,t)|\,dx\,dt
\nonumber
\\
&\qquad\le C\sum_{|\beta|=l}
\iint_{\widetilde CQ_{0,m}\times B^1}
\gamma(x,t)|D^\beta f(x,t)|\,dx\,dt.
\label{eq3.21}
\end{align}
}

Summing up estimate \eqref{eq3.21} over $m \in \mathbb{Z}^{n}$ and taking into account the finite multiplicity
of the cubes $\widetilde{C} Q_{k,m}^{n}$  (with fixed~$k$ and variable~$m$)   in view of~\eqref{eq3.20} we obtain

$$
\|f|W^{l}_{1}(\mathbb{R}^{n+1},\gamma)\| \le C(n,l,C_{\gamma},\Theta)\|\varphi|\widetilde{B}^{l}(\mathbb{R}^{n},\{\gamma_{k,m}\})\|.
$$

It remains to show that $\varphi=\operatorname{tr}|_{y=0} f$.
We fix an arbitrary cube~$Q$ in $\mathbb{R}^{n}$. Almost every point  $x\in\mathbb R^n$
is a~Lebesgue point of the function $\varphi$, because
$\varphi\in L^{\mathrm{loc}}_1(\mathbb R^n)$.
Hence, for almost all $x\in\mathbb R^n$,
$$
g_\delta(x):=\frac{1}{\delta^n}
\int_{x+\delta I}|\varphi(x')-\varphi(x)|\,dx'\to 0\quad
\text{as }\ \delta\to 0.
$$
Consequently, by the Lebesgue convergence theorem,
\begin{equation}
\int_{Q^n}|\varphi(x)-E_\delta[\varphi](x)|\,dx
\le C\int_{Q^n}g_{\widetilde{\delta}}(x)\,dx\to 0\quad
\text{as }\
\delta\to 0.
\label{eq3.22}
\end{equation}

Note that in the right-hand side of \eqref{eq3.22} we put $\widetilde{\delta}:=C(n,l,\Theta)\delta$.

From \eqref{eq3.22} and the definition of the function~$f$ it easily follows that $\varphi$ is the trace of the function $f$ on the hyperplane $x_{n+1}=0$.

The proof of Theorem~\ref{Th3.1}  is complete.

\begin{Remark} \rm
 In the case $l=1$, $\gamma(x_{1},..,x_{n+1})=|x_{n+1}|^{-\alpha}$, $\alpha \in (0,1)$, all the arguments employed in the proof of Theorem~3.1 remain valid.
 By $B^{l-1+\alpha}_{1,1}$ we shall denote the  classical Besov space. Using the fact that $\widetilde{B}^{l}(\mathbb{R}^{n},\{\gamma_{k,m}\})=B^{l-1+\alpha}_{1,1}$
 (the proof may be found in \cite{Tyu}, Remark 2.9) with $\gamma(x_{1},..,x_{n+1})=|x_{n+1}|^{-\alpha}$, $\alpha \in (\min\{1,l-1\},\max\{l,l-1\})$,
 we see that in this setting our result agrees with those obtained in Theorems 1.1, 1.2 of~\cite{Mi}. It is worth pointing out that in~\cite{Mi} it was assumed that the weight $\gamma=|x_{n+1}|^{-\alpha}$.
 Under this assumption the paper \cite{Mi} is capable of encompassing the cases  $\gamma \notin A^{\rm loc}_{1}(\mathbb{R}^{n+1})$.
 However, after some modifications of the proof of our Theorem 3.1 one may show that for $\alpha \in (\min\{1,l-1\},\max\{l,l-1\})$
 the trace of the space $W^{l}_{1}(\mathbb{R}^{n+1},|x_{n+1}|^{-\alpha})$ is the classical Besov space $B^{l-1+\alpha}_{1,1}$.
\end{Remark}

\section{The limiting case}

In this section we shall be concerned with the problem of complete description of the trace space of the Sobolev space $W^{l}_{1}(\mathbb{R}^{n+1}_{+},\gamma)$ with
$l=1$ and $\gamma \in A^{\rm loc}_{1}(\mathbb{R}^{n+1})$. We first note that this problem is equivalent to the problem of the description
of the trace space of the Sobolev space
$W^{1}_{1}(\mathbb{R}^{n+1},\gamma)$ on~$\mathbb{R}^{n}$. Indeed, this follows from the easily verified fact that the operator  of even extension from
$W^{1}_{1}(\mathbb{R}^{n+1}_{+},\gamma)$ into the space $W^{1}_{1}(\mathbb{R}^{n+1},\gamma)$ is continuous.

Before proceeding with precise statements, we first give a brief `heuristic' description of this problem in order to clarify, on the intuitive level, the principal impetuses for further
constructions.

Unfortunately, the Besov-type space of variable smoothness $\widetilde{B}^{1}(\mathbb{R}^{n},\{\gamma_{k,m}\})$ (considered in the previous section)
are poor candidate for the role of trace space if a~weight is only subject to the constraint $\gamma \in A^{\rm loc}_{1}(\mathbb{R}^{n+1})$.

It is not hard to see that  for $l=1$ estimate \eqref{eq3.9} fails in general, and hence one may not expect an estimate like \eqref{eq3.13}.
In addition to this technical impediment there are much deeper reasons for the unfitness of the spaces $\widetilde{B}^{1}(\mathbb{R}^{n},\{\gamma_{k,m}\})$.

Indeed, even in the case $\gamma \equiv 1$ the classical Gagliardo's result  shows that
$\operatorname{Tr}\left|_{t=0}W^{1}_{1}(\mathbb{R}^{n+1}_{+})\right.=L_{1}(\mathbb{R}^{n})$. In this case the space
$\widetilde{B}^{1}(\mathbb{R}^{n},\{\gamma_{k,m}\})$ coincides with the Besov space of smoothness zero
$B^{0,1}_{1,1}(\mathbb{R}^{n})$ (see \cite{Be}~for details). Next, Lemma~2 of~\cite{Be} implies, in particular, that
$\widetilde{B}^{1}(\mathbb{R}^{n},\{\gamma_{k,m}\}) \neq L_{1}(\mathbb{R}^{n})$ for $\gamma \equiv 1$.
So, the trace space contains functions with inappropriate smoothness properties. It is also worth pointing out that, according to Peetre~\cite{Pe},
the extension operator   $\operatorname{Ext}:L_{1}(\mathbb{R}^{n}) \to W^{1}_{1}(\mathbb{R}^{n+1}_{+})$
(which is the right inverse of the trace operator) cannot be linear.

On the other hand,
for the weight $\gamma=\gamma(x,t)=|t|^{-\varepsilon}$,  $\varepsilon \in (0,1)$, the methods of the previous section also work!
To this aim one needs to slightly modify estimate \eqref{eq3.9}. In spite of the fact that  $l-1=0$, we succeed in achieving a~`geometric rate' on account of the
fact that $\inf\limits_{t \in (0,2^{-k})}t^{-\varepsilon} \geq 2^{\varepsilon}\inf\limits_{t \in (0,2^{-k+1})}t^{-\varepsilon}$.

As a good candidate for the trace space in the general case one should consider a~space whose elements are able to  appreciably change their
smoothness characteristics when transiting from a~point to a~point, because the `rate of decay' of a~weight may be substantially different at different points.
As distinct from the case  $l > 1$, in which, roughly speaking, the trace space is `quasi-homogeneous', in the case $l=1$ the trace space turns out to be `essentially nonhomogeneous'.

In the case $l > 1$ a~sufficiently rapidly growing geometric progression $\{2^{kl}\}$
helped to control the strong inhomogeneity of a~weight.
However, the limiting case $p=l=1$ calls for a~more subtle analysis of the local behaviour of the weight near each point of the hyperplane on which the trace is considered.

An important step in this analysis is the construction of a~special system of tilings of the space $\mathbb{R}^{n}$.
This system  of tilings will replace the standard system of tilings of the space $\mathbb{R}^{n}$ composed of all dyadic cubes
numbered by indexes  $(k,m) \in \mathbb{Z}_{+} \times \mathbb{Z}^{n}$. The cubes in our special system of tilings will be numbered by indexes
$(k,m) \in A \subset \mathbb{Z}_{+} \times \mathbb{Z}^{n}$. Here, the algorithm for construction of the index set~$A$ is based on combinatorial arguments and is nonlinear.
Namely, the set~$A$ depends not only on the weight~$\gamma$, but also on the function~$f$.

In this section we shall denote by $Q$ (respectively, $\overline{Q}$) an open (closed) cube in the space~$\mathbb{R}^{n}$.

 \begin{Def} \rm
 \label{Def4.1}
Assume that we are given a set of dyadic closed cubes $T=\{\overline{Q}_{\alpha}\}_{\alpha \in A}$, $A \subset \mathbb{Z}_{+} \times \mathbb{Z}^{n}$,
 in which different cubes have disjoint interiors and $\mathbb{R}^{n}=\bigcup\limits_{\alpha \in A}\overline{Q}_{\alpha}$. We shall call this family a~\textit{tiling  of the space}~$\mathbb{R}^{n}$.
\end{Def}

\begin{Def}\rm
 \label{Def4.2}
A~tiling $T'=\{\overline{Q}_{\alpha}\}_{\alpha \in A'}$ will be said to succeed a tiling $T=\{\overline{Q}_{\alpha}\}_{\alpha \in A}$ (written $T' \succ T$)
if each cube $\overline{Q}_{\alpha'}$, $\alpha' \in A'$, of the tiling~$T'$ is contained in some cube $\overline{Q}_{\alpha}$, $\alpha \in A$, of the tiling~$T$.
\end{Def}

\begin{Def} \rm
\label{Def4.3}
Assume that for any $s \in \mathbb{Z}_{+}$ we have a~tiling $T^{s}=\{\overline{Q}^{s}_{\alpha}\}_{\alpha \in A^{s}}$, $A^{s} \subset \mathbb{Z}_{+} \times \mathbb{Z}^{n}$
of the space~$\mathbb{R}^{n}$. Assume also that $T^{s+1} \succ T^{s}$ for $s \in \mathbb{Z}_{+}$. Then the set $T=\{T^{s}\}=\{T^{s}\}^{\infty}_{s=0}$
will be called a~system of tilings of the space~$\mathbb{R}^{n}$.
\end{Def}

The next lemma is an important combinatorial instrument required in the definition of the trace space.

\begin{Lm}
\label{Lm4.1}
Let $\{\overline{Q}_{\alpha}\}_{\alpha \in A}$ be a tiling of the space $\mathbb{R}^{n}$. Then, for each number $\lambda = 2^{-k_{0}}$, $k_{0} \in \mathbb{Z}_{+}$,  there exists
an index set $\widetilde{A} \subset A$ such that
\begin {list}{}{\itemsep=0pt\topsep=2pt\parsep=0pt}
\item[\rm  1)] $\mathbb{R}^{n} = \bigcup\limits_{\alpha \in \widetilde{A}}\widetilde{Q}_{\alpha}$, $\widetilde{Q}_{\alpha}:=(1+\lambda)Q_{\alpha}$,
\item[\rm  2)] any point $x \in \mathbb{R}^{n}$ lies in at most  $(n+1)2^{n}$ cubes from the family $\{\widetilde{Q}_{\alpha}\}_{\alpha \in \widetilde{A}}$,
\item[\rm  3)] if $\widetilde{Q}_{\alpha} \bigcap \widetilde{Q}_{\alpha'} \neq \varnothing$, then
$|\widetilde{Q}_{\alpha} \bigcap \widetilde{Q}_{\alpha'}| \geq C (n,\lambda) \min\{|Q_{\alpha'}|,|Q_{\alpha}|\}$,
\item[\rm  4)] every cube $\widetilde{Q}_{\alpha}, \alpha \in \widetilde{A}$ is not contained in $\bigcup\limits_{\alpha' \in \widetilde{A}, \alpha' \neq \alpha}\widetilde{Q}_{\alpha'}$.
\end{list}
\end{Lm}

\textbf{Proof.}
Since the set $A$ is countable, we can enumerate all the cubes $\{\overline{Q}_{\alpha}\}_{\alpha \in A}$ by natural numbers: $\{\overline{Q}_{\alpha}\}_{\alpha \in A}=\{\overline{Q}_{\alpha_{i}}\}_{i=1}^{\infty}$. We set $S^{0}:=A$. Let $\alpha_{i_{1}}$ be the first index for which  $\widetilde{Q}_{\alpha_{i_{1}}} \subset \bigcup
\limits_{i \in \mathbb{N}, i\neq i_{1}}\widetilde{Q}_{\alpha_{i}}$. If there is no such index, then we set  $\widetilde{A}=A$ and complete the construction.
We exclude the cube  $Q_{\alpha_{1}}$ from our system and consider the index set $S^{1}:=A \setminus \{\alpha_{i_{1}}\}$.  It is clear that
$\mathbb{R}^{n}=\bigcup\limits_{\alpha \in S^{1}}\widetilde{Q}_{\alpha}$. Assume that we have already constructed indexes  $i_{1} < ... < i_{k}$ and sets
$S^{1} \supset .. \supset S^{k}$. Let  $i_{k+1} > i_{k}$ be the first natural number for which  $\widetilde{Q}_{\alpha_{i_{k+1}}} \subset \bigcup\limits_{\alpha \in S^{k}, \alpha \neq \alpha_{i_{k+1}}}\widetilde{Q}_{\alpha}$. If there is no such number, then we put $\widetilde{A}=S^{k}$ and complete the construction.
We exclude the cube  $\widetilde{Q}_{\alpha_{i_{k+1}}}$ from our system and consider the index set $S^{k+1}:=S^{k} \setminus \{\alpha_{i_{k+1}}\}$.
It is easily checked that $\mathbb{R}^{n}=\bigcup\limits_{\alpha \in S^{k+1}}\widetilde{Q}_{\alpha}$. This being so,
either the set $\widetilde{A}$ will be obtained in a~finite number of steps or we get an increasing sequence of natural numbers $\{i_{k}\}_{k=1}^{\infty}$ and
a~sequence of sets $S^{1}\supset \dots \supset S^{k} \supset \dots$. Let  $A=\bigcap\limits_{k=1}^{\infty}S^{k}$. We claim that $\widetilde{A}$ is the required index set.

By the construction,  $\mathbb{R}^{n}=\bigcup\limits_{\alpha \in S^{k}}\widetilde{Q}^{s}_{\alpha}$ for each $k \in \mathbb{N}$, and hence
$\mathbb{R}^{n}=\bigcup\limits_{\alpha \in \widetilde{A}}\widetilde{Q}^{s}_{\alpha}$. This proves assertion~1).

Let us prove assertion 4). Assume there is a~cube  $\widetilde{Q}_{\alpha_{i_{k_{0}}}}$ such that  $\widetilde{Q}_{\alpha_{i_{k_{0}}}} \subset \bigcup\limits_{\alpha \in \widetilde{A} \setminus \{\alpha_{i_{k_{0}}}\}}\widetilde{Q}_{\alpha}$.
But then there exists a biggest number $0 \le k'_{0} \le k_{0}$ and a~set $S^{k'_{0}}$ for which
$\widetilde{Q}_{\alpha_{i_{k_{0}}}} \subset \bigcup\limits_{\alpha \in S^{k'_{0}}}\widetilde{Q}_{\alpha}$, contradicting the construction of $\widetilde{A}$.

2) We claim that the overlapping multiplicity is at most $(n+1)2^{n}$.
Indeed, we fix an arbitrary point $x_{0} \in \mathbb{R}^{n}$ and estimate the number of cubes from the family
$\{\widetilde{Q}_{\alpha}\}_{\alpha \in \widetilde{A}}$ that contain this point. In doing so we shall modify one trick from Lemma~1.1 of~\cite{Guz}.
Namely, we draw through the point $x_{0}$ the planes that are parallel to the coordinate planes. This will give us $2^{n}$ quadrants (closed!) with vertex at $x_{0}$.
We fix arbitrary quadrant and consider the cubes that contain~$x_{0}$ and whose centers lie in this quadrant.
Clearly, the lemma will be proved once we show that there are at most $(n+1)$ such cubes.

Assume the contrary. Note that if the centres of the cubes $\widetilde{Q}_{\alpha} \ni x_{0}$ and $\widetilde{Q}_{\alpha'} \ni x_{0}$ lie in the same quadrant,
then the center of one cube lies in the other cube. This implies, in particular, that either
$Q_{\alpha} \subset \widetilde{Q}_{\alpha'}$ or $Q_{\alpha'} \subset \widetilde{Q}_{\alpha}$ (inasmuch as $\lambda = 2^{-k}$, $k \in \mathbb{Z}_{+}$).
It follows that if $r(Q_{\alpha})=r(Q_{\alpha'})$ then these two cubes coincide. Hence, we may assume that the side lengths of the cubes containing the point $x_{0}$
and whose centres lie in the same quadrant, are strictly decreasing. Then, numbering these cubes in decreasing size, we see that the centre of the next cube (in the order of decreasing size) is contained in its direct predecessor.
As a~result, the centre of the cube with number $n+2$ (which we denote by $\widetilde{Q}_{\alpha_{0}}$) will be contained in at least~$n+1$  cubes
from the family  $\{\widetilde{Q}_{\alpha}\}_{\alpha \in \widetilde{A}}$. We claim that such a~case is never realized (we shall obtain a~contradiction with the
algorithm for choosing the cubes).

The key observation here is that  $\overline{Q}_{\alpha_{0}} \subset \widetilde{Q}_{\alpha'}$
if and only if $\widetilde{Q}_{\alpha_{0}} \subset \widetilde{Q}_{\alpha'}$. Hence if $Q_{\alpha_{0}} \subset \widetilde{Q}_{\alpha'}$ (for $\alpha_{0},\alpha' \in \widetilde{A}$),
then the closed cube $\overline{Q}_{\alpha_{0}}$ cannot wholly lie in the cube $\widetilde{Q}_{\alpha'}$ (because otherwise $\widetilde{Q}_{\alpha_{0}} \subset \widetilde{Q}_{\alpha'}$ and the cube $\widetilde{Q}_{\alpha_{0}}$ will be excluded during the construction of the set $\widetilde{A}$).

So, having a fixed quadrant and a~cube $\widetilde{Q}_{\alpha_{0}}$ with center in this quadrant, we estimate the number of cubes
$\widetilde{Q}_{\alpha}$, $\alpha \in \widetilde{A}$, whose centers lie inside this quadrant, which contains the cube $Q_{\alpha_{0}}$, but which
do not contain the cube $\overline{Q}_{\alpha_{0}}$.
We will show that there are at most~$n$ such cubes.

We note that the facets of any~$Q$ can be canonically labeled by natural numbers from 1 to~$2n$. Assume that
a~cube $\widetilde{Q}_{\alpha'} \supset Q_{\alpha_{0}}$ does not contain the cube $\overline{Q}_{\alpha_{0}}$ and $r(Q_{\alpha'}) > r(Q_{\alpha_{0}})$.
Then, for some $i \in \{1,\dotsc,2n\}$, the $i$th facet of the cube $\overline{Q}_{\alpha_{0}}$ lies in the $i$th facet of the cube $\overline{\widetilde{Q}}_{\alpha'}$,
for otherwise we would get the inclusion  $\overline{Q}_{\alpha_{0}} \subset \widetilde{Q}_{\alpha'}$
(because $Q_{\alpha_{0}} \subset \widetilde{Q}_{\alpha'}$), which contradicts the construction. Next, if the $i$th facet of the cube $\overline{\widetilde{Q}}_{\alpha'}$ contains the $i$th
 facet of the cube $\overline{Q}_{\alpha_{0}}$,
 then there is no other cube $\widetilde{Q}_{\alpha''}$ containing the cube $Q_{\alpha_{0}}$ (whose center lies in the quadrant under consideration!)
 which has such a~property.

Indeed, let $Q_{\alpha_{0}} \subset \widetilde{Q}_{\alpha''}$, $Q_{\alpha_{0}} \subset \widetilde{Q}_{\alpha'}$ and
$r(Q_{\alpha_{0}}) < r(Q_{\alpha''}) < r(Q_{\alpha'})$. We claim that if the $i$th facet of the cube $\overline{Q}_{\alpha_{0}}$ is contained in the $i$th
 facet of the cube  $\overline{\widetilde{Q}}_{\alpha'}$, then the $i$th facet of the cube $\overline{Q}_{\alpha_{0}}$ is not contained in the  $i$th facet of the cube
 $\overline{\widetilde{Q}}_{\alpha''}$.

 From the construction of the index set $\widetilde{A}$ it follows that the closed cube $\overline{Q}_{\alpha_{0}}$
 is not wholly contained in the cube $\widetilde{Q}_{\alpha'}$ and in the cube $\widetilde{Q}_{\alpha''}$.
 Hence,  $Q_{\alpha_{0}}\bigcap Q_{\alpha'} = \emptyset$, $Q_{\alpha_{0}}\bigcap Q_{\alpha''} = \emptyset$.

Consider dyadic cubes with side length $\frac{\lambda}{2} r(Q_{\alpha'})$ lying in the set $\widetilde{Q}_{\alpha'} \setminus \overline{Q}_{\alpha'}$.
Since $Q_{\alpha_{0}} \subset \widetilde{Q}_{\alpha'}$ and since $Q_{\alpha_{0}} \bigcap Q_{\alpha'} = \varnothing$, among the above cubes there exists a~unique
dyadic cube $Q_{\beta_{0}} \supset Q_{\alpha_{0}}$ (note that $r(Q_{\beta_{0}}) =\frac{\lambda}{2} r ({Q}_{\alpha'})$).
Besides, the $i$th facet of the cube $\overline{Q}_{\alpha_{0}}$ is contained in the $i$th facet of the cube
$\overline{Q}_{\beta_{0}}$ (because by the assumption the $i$th facet of the cube $\overline{Q}_{\alpha_{0}}$ lies in the $i$th facet of the cube
$\overline{\widetilde{Q}}_{\alpha'}$).

For further purposes we shall need the following key observation. Assume that we are given two arbitrary cubes  $Q_{\alpha'}$ and $Q_{\alpha''}$.
Then the distance between the hyperplanes containing the $i$th facets of these cubes is either zero or is not smaller than the side length of the smallest of these
2~cubes.

Now we consider two cases. In the first case  $r(Q_{\beta_{0}}) \le r(Q_{\alpha''})$ and $Q_{\alpha_{0}} \subset \widetilde{Q}_{\alpha''}$.
If the $i$th facet of the cube $\overline{Q}_{\alpha'}$ and the $i$th facet of the cube $\overline{Q}_{\alpha''}$ lie in the same hyperplane, then it is clear that the $i$th facet of the cube
$\overline{Q}_{\alpha_{0}}$ cannot simultaneously lie in the same hyperplane with the $i$th facet of the cube $\overline{\widetilde{Q}}_{\alpha'}$ and in the same hyperplane  with the $i$th
facet of the cube $\overline{\widetilde{Q}}_{\alpha''}$ (as required).
If, however, the $i$th facet of the cube $\overline{Q}_{\alpha'}$ and the $i$th facet of the cube $\overline{Q}_{\alpha''}$ do not lie in the same hyperplane, then the distance
between the hyperplanes that contain these facets is not smaller than $r(Q_{\alpha''}) \geq r(Q_{\beta_{0}})$. But in this case
it follows by simple geometrical considerations that the distance between the hyperplanes  containing the $i$th facet of the cube $\overline{\widetilde{Q}}_{\alpha'}$and the  $i$th facet of the
cube $\overline{\widetilde{Q}}_{\alpha''}$ is positive. Hence, the $i$th facets of the cubes  $\widetilde{Q}_{\alpha_{0}}$, $\overline{\widetilde{Q}}_{\alpha'}$, $\overline{\widetilde{Q}}_{\alpha''}$
do not lie in the same hyperplane.

In the second case  $r(Q_{\beta_{0}}) > r(Q_{\alpha''})$. If the distance between the hyperplanes containing the $i$th facets of the
cubes $\overline{Q}_{\beta_{0}}$ and $\overline{Q}_{\alpha''}$ is positive, then by the above observation and since $\frac{\lambda}{2} \le \frac{1}{2}$, it follows that the distance
 from the hyperplane containing the $i$th facet of the cube $\overline{\widetilde{Q}}_{\alpha''}$ to the hyperplane containing the $i$th facet of the cube $\overline{Q}_{\beta_{0}}$ (and hence $\overline{\widetilde{Q}}_{\alpha'}$) is positive.
If now the $i$th facets of the cubes $\overline{Q}_{\beta_{0}}$ and $\overline{Q}_{\alpha''}$ lie on one hyperplane, then the $i$th facet of the cube $\overline{\widetilde{Q}}_{\alpha''}$
does not lie in the same hyperplane with them.

In all cases considered above, we see that the $i$th ~facet of the cube $\overline{Q}_{\alpha_{0}}$ is contained in the $i$th facet of the cube  $\overline{Q}_{\beta_{0}}$
(and hence,  $\widetilde{Q}_{\alpha'}$), but is not contained in the $i$th facet of the cube $\overline{\widetilde{Q}}_{\alpha''}$.

Let $j=j(i)$ be the index corresponding to the facet which is parallel to the $i$th facet (recall that the facet are labeled in the canonical way and that the
labeling is the same for each cube). It now remains to note that if the $i$th facet of the cube $\overline{Q}_{\alpha_{0}}$ lies in the $i$th facet of the cube
$\overline{\widetilde{Q}}_{\alpha'} \supset Q_{\alpha_{0}}$, the $j$th facet of the cube  $\overline{Q}_{\alpha_{0}}$ lies in the $j$th facet of some
$\overline{\widetilde{Q}}_{\alpha''} \supset Q_{\alpha_{0}}$, and besides, $r(Q_{\alpha'}),r(Q_{\alpha''}) > r(Q_{\alpha_{0}})$, then the centers of the
cubes $\widetilde{Q}_{\alpha'}$ and $\widetilde{Q}_{\alpha''}$ cannot lie in the same quadrant.

Let us now prove assertion 3). Assume that $\widetilde{Q}_{\alpha} \bigcap \widetilde{Q}_{\alpha'} \neq \emptyset$. We set
$l_{0}=\frac{\lambda}{2} \min\{r(Q_{\alpha}),r(Q_{\alpha'})\}$. Then the cube $\widetilde{Q}_{\alpha}$ and the  $\widetilde{Q}_{\alpha'}$ can be represented as a~union of
dyadic cubes (possibly containing portions of their boundaries) with side length $l_{0}$. But open dyadic cubes of the same size length are either disjoint or equal. Hence,
there exists at least one cube with side length $l_{0}$ which is contained both in the cube $\widetilde{Q}_{\alpha}$ and in the cube $\widetilde{Q}_{\alpha'}$
(because the intersection of such cubes is nonempty). Now the required estimate $|\widetilde{Q}_{\alpha} \bigcap \widetilde{Q}_{\alpha'}| \geq (l_{0})^{n} \geq C(n,\lambda)\min\{|Q_{\alpha}|,|Q_{\alpha'}|\}$
is clear.

\textbf{Notations}. We shall frequently use the following notation. Given a fixed parameter $\lambda=2^{-k}$, $k\in \mathbb{Z}_{+}$, and a cube~$Q$ in $\mathbb R^n$, we set $\widetilde Q:=(1+\lambda)Q$. Given $k \in \mathbb{Z}_{+}$, $m \in \mathbb{Z}^{n}$, we set
\begin{gather*}
\Pi_{k,m}:=Q_{k,m}\times(0,r(Q_{k,m})), \qquad \widetilde{\Pi}_{k,m}:=\widetilde{Q}_{k,m}\times(0,r(Q_{k,m})),
\\
\widehat{\gamma}_{k,m}:=\gamma_{\Pi_{k,m}}, \qquad   \gamma_{k,m}:= (r(Q^{s}_{\alpha}))^{n+1}\widehat{\gamma}_{k,m}.
\end{gather*}

If $T=\{Q_\alpha\}_{\alpha\in A}$ is a~tiling of~$\mathbb R^n$, then by $\tilde A$~we shall denote the index set which was constructed in Lemma~\ref{Lm4.1}.

Assume we are given a system of tilings $T=\{T^{s}\}$ of the space $\mathbb{R}^{n}$ and a fixed parameter $\lambda=2^{-k}$, $k\in \mathbb{Z}_{+}$.
Applying Lemma~\ref{Lm4.1} for each $s \in \mathbb{Z}_{+}$ to the tiling $T^{s}$, we obtain the covering
$\Xi^{s}$ of the space $\mathbb{R}^{n}$ by cubes  $\{\widetilde{Q}^{s}_{\alpha}\}_{\alpha \in \widetilde{ A^{s}}}$.

In the cases when we know that $\alpha=(k,m) \in A^{s} \subset \mathbb{Z}_{+} \times \mathbb{Z}^{n}$, then instead of $\widehat{\gamma}_{k,m}$, $\gamma_{k,m}$, $\widetilde{\Pi}_{k,m}$, $\Pi_{k,m}$,
we shall write, respectively, $\widehat{\gamma}^{s}_{\alpha}$, $\gamma^{s}_{\alpha}$, $\widetilde{\Pi}^{s}_{\alpha}$, $\Pi^{s}_{\alpha}$.

For a~function $\varphi \in L^{\rm loc}_{1}(\mathbb{R}^{n})$ and a~given system tilings $T=\{T^{s}\}$, we denote
$$
\varphi^{s}_{\alpha} = \frac{1}{|\widetilde{Q}^{s}_{\alpha}|}\int\limits_{\widetilde{Q}^{s}_{\alpha}}\varphi(x)\,dx, \qquad s \in \mathbb{Z}_{+}, \ \ \alpha \in \widetilde{A}^{s}.
$$

By $\widetilde{q}$ we shall denote the smallest of $C \geq 1$ for which
$$
\frac{1}{8|\Pi_{k,m}|}\int\limits_{8 Q_{k,m}}\int\limits_{0}^{r(Q_{k,m})}\gamma(x,t)\,dt\,tx\le
C\widehat{\gamma}_{k,m'},
$$
where $k\in \mathbb Z_+$, $m\in \mathbb Z^n$, and $Q_{k,m'}\subset 8Q_{k,m}$.

Let $q:=16 \widetilde{q} C_{\gamma} 2^{n+1} $. From the definition of $\widetilde{q}$ and \eqref{eq2.2} we have

\begin{equation}
\label{eq.q}
\widehat{\gamma}_{k,m} \le \frac{q}{2} \widehat{\gamma}_{k,m'}, \quad k\in \mathbb Z_+, |m_{i}-m'_{i}| \le 1, i \in \{1,..,n\},
\end{equation}

\begin{equation}
\label{eq.q'}
\widehat{\gamma}_{k,m} \le \frac{q}{2}\widehat{\gamma}_{k+1,m'} \le (\frac{q}{2})^{2} \widehat{\gamma}_{k,m}, \quad k\in \mathbb Z_+, m,m' \in \mathbb{Z}^{n}, Q_{k+1,m'} \subset Q_{k,m}.
\end{equation}

The role of the parameter $q$ will be transparent at Step~1 of the proof of Theorem 4.1.

\begin{Def} \rm
\label{Def4.4}
Let $c_{1},c_{2} \geq 1$. A~ system of tilings $T=\{T^{s}\}=\{T^{s}\}_{s=0}^{\infty}(c_{1},c_{2})$ of the space~$\mathbb{R}^{n}$  ($T^{s}=\{\overline{Q}^{s}_{\alpha}\}_{\alpha \in A^{s}}$,
where $s \in \mathbb{Z}_{+}$), is called \textit{admissible for a~weight}~$\gamma$ if, for each $s \in \mathbb{Z}_{+}$, the following conditions are satisfied:

1) if  $\widetilde{Q}^{s}_{\alpha}\bigcap\widetilde{Q}^{s}_{\alpha'} \neq \varnothing$, then
$\widehat{\gamma}^{s}_{\alpha} \le c_{1} \widehat{\gamma}^{s}_{\alpha'}$  for $\alpha,\alpha' \in A^{s}$;

2) if $Q^{s+1}_{\alpha'} \subset Q^{s}_{\alpha}$, then $\widehat{\gamma}^{s}_{\alpha} \le c_{2}\widehat{\gamma}^{s+1}_{\alpha'}$ and $\widehat{\gamma}^{s+1}_{\alpha'} \le c_{2} \widehat{\gamma}^{s}_{\alpha}$ for $\alpha \in A^{s}$, $\alpha' \in A^{s+1}$;

3) $\max\{r(Q^{s+1}_{\alpha'}): x \in \widetilde{Q}^{s+1}_{\alpha'}\} \le \frac{1}{2} \min\{r(Q^{s}_{\alpha}): x \in \widetilde{Q}^{s}_{\alpha}\}$
for any $x \in \mathbb{R}^{n}$;

4) $r(Q^{s}_{\alpha}) \geq 2^{-l_{s}}$, $\alpha \in A^{s}$, for some strictly increasing sequence of nonnegative integer numbers $\{l_{j}\}_{j=0}^{\infty}$ with $l_{0}=0$

\end{Def}

\begin{Th}
\label{Th4.1}
Let a weight
$\gamma \in A^{\rm loc}_{1}(\mathbb{R}^{n+1})$, $\lambda=1+2^{-k}$ $(k \in \mathbb{Z}_{+})$. Then there exist constants $c_{1}(n,C_{\gamma}),c_{2}(n,C_{\gamma}) \geq 1$ such that, for any function
$f \in W^{1}_{1}(\mathbb{R}^{n+1}_{+},\gamma)$, there exists a~{system of tilings} $T=\{T^{s}\}(c_{1},c_{2})$ of the space~$\mathbb{R}^{n}$
that is admissible for the weight~$\gamma$ and is such that
\begin{equation}
\begin{split}
\label{eq4.1}
&\sum\limits_{m \in \mathbb{Z}^{n}}\widehat{\gamma}_{0,m}\|\varphi|L_{1}(Q_{0,m})\| + \sum\limits_{s=1}^{\infty}\sum\limits_{\alpha \in \widetilde{A}^{s}} \widehat{\gamma}^{s}_{\alpha}\int\limits_{\widetilde{Q}^{s}_{\alpha}}|\varphi^{s}_{\alpha}-\varphi(x)|\,dx \le C(n,C_{\gamma},\lambda,c_{1},c_{2})\|f|W^{1}_{1}(\mathbb{R}^{n+1}_{+},\gamma)\|.
\end{split}
\end{equation}
\end{Th}

\textbf{Proof.}
\textit{Step $1$.} For each function $f \in W^{1}_{1}(\mathbb{R}^{n+1}_{+},\gamma)$ we construct a~required system  of tilings of the space~$\mathbb{R}^{n}$
that is admissible for the weight~$\gamma$.

First we construct an auxiliary system of tilings $\{\mathring T^{s}\}$ that will satisfy only properties 1), 2) and~4) of Definition~\ref{Def4.4}. Next, for $r \in \mathbb{N}$, $r > 1$
we choose a~required subsystem  $\{T^{s}\}:=\{\mathring T^{rs}\}$ of the system  $\{\mathring T^{s}\}$.

Let $\{l_{j}\}_{j=1}^{\infty}$ be a strictly increasing sequence of nonnegative integer numbers such that $l_{0}=0$ and
\begin{equation}
\label{eq4.2}
\|f|W^{1}_{1}(\mathbb{R}^{n}\times (0,2^{-l_{j+1}}))\| \le \frac{1}{2} \|f|W^{1}_{1}(\mathbb{R}^{n}\times (0,2^{-l_{j}}))\|, \qquad j \in \mathbb{Z}_{+}.
\end{equation}

We construct the required system of tilings  $\{\mathring T^{k}\}$ by induction.

\textit{Induction basis.} We first build a tiling $\mathring T^{0}$. To do so we put $\mathring T^{0}:=\{\overline{Q}_{0,m}\}_{m \in \mathbb{Z}^{n}}$,
$\mathring A^{0}:=\{0\} \times \mathbb{Z}^{n}$, and for each $m \in \mathbb{Z}^{n}$ we paint the cube $\overline{Q}_{0,m}$ yellow.

\textit{Induction step.} Assume that for $s \in \mathbb{Z}_{+}$ the tiling $\mathring T^{s}=\{\overline{Q}^{s}_{\alpha}\}_{\alpha \in \mathring A^{s}}$ is constructed.
Let us construct the tiling $\mathring T^{s+1}$.
We fix a~cube $\overline{Q}^{s}_{\alpha}$ for $\alpha \in \mathring A^{s}$. Suppose that $\widehat{\gamma}^{s}_{\alpha} \in [q^{j},q^{j+1})$
for some $j \in \mathbb{Z}$. We decompose the cube $\overline{Q}^{s}_{\alpha}$ into dyadic cubes
($\overline{Q}_{k,m}$, say) of twice smaller size. Among these cubes, we select those satisfying the estimate $\widehat{\gamma}_{k,m} > q^{j+1}$ and paint them blue.
Note that in view of estimate \eqref{eq.q'} we have $\widehat{\gamma}_{k,m} \in (q^{j+1},\frac{q^{j+2}}{2}]$ (it is important here that the parameter $q$ is sufficiently large!).
We decompose the remaining cubes into the cubes $\overline{Q}_{k+1,m'}$, select  those for which $\widehat{\gamma}_{k+1,m'} > q^{j+1}$ and paint these cubes blue.
This process is repeated until the side length of a~cube will be $2^{-l_{s+1}}$. In this case we either have a~tiling of the cube $\overline{Q}^{s}_{\alpha}$ consisting
of only blue cubes or there will be cubes  $\overline{Q}_{l_{s+1},m''} \subset \overline{Q}^{s}_{\alpha}$ for which  $\widehat{\gamma}_{l_{s+1},m''} \le q^{j+1}$.
In the latter case, we paint these cubes $\overline{Q}_{l_{s+1},m''}$ yellow. The resulting tiling of the cube $\overline{Q}^{s}_{\alpha}$ will be composed of the so-chosen blue cubes and
the remaining yellow cubes. Combining the corresponding tilings of the cubes $\overline{Q}^{s}_{\alpha}$ over all $\alpha \in \mathring A^{s}$, we obtain the
tiling  $\mathring T^{s+1}$ of the space~$\mathbb{R}^{n}$. By  $\mathring A^{s+1}$ we shall denote the set of pairs of indices $(k,m) \in \mathbb{Z}_{+} \times \mathbb{Z}^{n}$ for which
$\overline{Q}_{k,m} \in \mathring T^{s+1}$.

Clearly, for each $s \in \mathbb{Z}_{+}$, the tiling  $\mathring T^{s}$ is composed of at most countable set of dyadic cubes.

If we apply Lemma 4.1 for each $s \in \mathbb{N}$ to the tiling  $\mathring T^{s}$, we obtain a~covering
$\mathring \Xi^{s}$ of the space~$\mathbb{R}^{n}$ by cubes $\{\widetilde{Q}^{s}_{\alpha}\}_{\alpha \in \widetilde{\mathring A^{s}}}$.

We next check that the system of tilings $\{\mathring T^{s}\}$ satisfies conditions 1), 2) and~4) of Definition \ref{Def4.4} (in which the index sets $A^{s}$
should be replaced by $\mathring A^{s}$).

Condition 4) is easily seen to hold.

We claim that Condition 1) of Definition \ref{Def4.4} is satisfied with constant  $c_{1}=q^{3}$. Let
$\widetilde{Q}^{s}_{\alpha} \bigcap \widetilde{Q}^{s}_{\alpha'} \neq \varnothing$
for $\alpha,\alpha' \in \mathring A^{s}$. Assume that $\widehat{\gamma}^{s}_{\alpha'} > q^{3}\widehat{\gamma}^{s}_{\alpha}$. For any cube $\overline{Q}^{s}_{\alpha}$, we let $b(\overline{Q}^{s}_{\alpha})$ denote the number of blue cubes $\overline{Q}^{j}_{\alpha'} \supset \overline{Q}^{s}_{\alpha}$ (for $j \le s$ and $\alpha' \in \mathring A^{j}$).
From our assumption it follows that there exists a~natural $k_{0} > 1$ such that the number of blue cubes containing the cube $\overline{Q}^{s}_{\alpha'}$ is greater by $k_{0}$ than the
number of blue cubes containing the cube $\overline{Q}^{s}_{\alpha}$. But then there exist a~blue cube $\overline{Q}^{s_{0}}_{\alpha'_{0}} \supset Q^{s}_{\alpha'}$ and a~yellow cube $\overline{Q}^{s_{0}}_{\alpha_{0}} \supset Q^{s}_{\alpha}$ such that
$\widehat{\gamma}^{s_{0}}_{\alpha'_{0}} \geq q^{k_{0}-1} \widehat{\gamma}^{s_{0}}_{\alpha_{0}}$. By the construction,
$r(Q^{s_{0}}_{\alpha_{0}}) \le r(Q^{s_{0}}_{\alpha'_{0}})$. Besides, $\widetilde{Q}^{s_{0}}_{\alpha_{0}}\bigcap\widetilde{Q}^{s_{0}}_{\alpha'_{0}} \neq \varnothing$
by  $\widetilde{Q}^{s}_{\alpha}\bigcap\widetilde{Q}^{s}_{\alpha'} \neq \varnothing$.
It follows that  $Q^{s_{0}}_{\alpha_{0}} \subset 8 Q^{s_{0}}_{\alpha'_{0}}$ (because $\lambda \le 1$) , and hence, $\widehat{\gamma}^{s_{0}}_{\alpha_{0}} \geq \frac{2}{q}\widehat{\gamma}^{s_{0}}_{\alpha'_{0}}$. A~contradiction is reached.

We now check condition~2). Let $Q^{s}_{\alpha}$ be the parent of the cube $Q^{s+1}_{\alpha'}$. By the construction of the system of  tilings, we have
$\widehat{\gamma}^{s}_{\alpha} \le  \widehat{\gamma}^{s+1}_{\alpha'}$ and $\widehat{\gamma}^{s+1}_{\alpha'} \le  q\widehat{\gamma}^{s}_{\alpha}$.

Let $r \in \mathbb{N}$, $r \geq 5$. Consider the system of tilings $\{T^{s}\}:=\{\mathring T^{rs}\}$ and define $A^{s}:=\mathring A^{rs}$. Clearly,
the system of tilings $\{T^{s}\}$ satisfies conditions 1), 2) (with the constants $c_{1}=q^{3}$, $c_{2}=q^{r}$) and~4) of Definition~\ref{Def4.4}.

Let us check condition~3) of Definition \ref{Def4.4}.  To this aim we fix a~point $x \in \mathbb{R}^{n}$. Let $\widetilde{Q}^{r(s+1)}_{\alpha} \ni x$ be a~cube
with largest side length among the set of cubes $\{\widetilde{Q}^{r(s+1)}_{\alpha}\}_{\alpha \in A^{r(s+1)}}$ that contain the point~$x$ (this cube may not be unique).
Let  $\widetilde{Q}^{rs}_{\alpha'} \ni x$ be a~cube of smallest side length among all cubes from the family $\{\widetilde{Q}^{rs}_{\alpha}\}_{\alpha \in A^{rs}}$,
of which each contains the point~$x$ (the cube $\widetilde{Q}^{rs}_{\alpha'} \ni x$ may also be not unique). Consider the following chain of nested dyadic
cubes $\overline{Q}^{r(s+1)}_{\alpha} \subset \dots. \subset \overline{Q}^{rs}_{\alpha''}$ (in this chain each succeeding dyadic cube is a~unique parent of its predecessor).
If this chain contains at least one yellow cube, then we have the result required. Indeed, by the construction, for any  $Q^{rs}_{\alpha}$, $\alpha \in \widetilde{A}^{s}$ (and hence, for $Q^{rs}_{\alpha'}$)
we have the estimate  $r(Q^{rs}_{\alpha}) \geq 2^{-l_{rs}}$. If the cube $\overline{Q}^{r(s+1)}_{\alpha}$ is yellow, then
$r(\overline{Q}^{r(s+1)}_{\alpha})=2^{-l_{r(s+1)}} < \frac{1}{2} 2^{-l_{rs}}$.
If another cube of the above chain is yellow, then the side length of this cube is clearly smaller or equal than $2^{-l_{rs}}$. The cube $\overline{Q}^{r(s+1)}_{\alpha}$
lying strictly inside it and hence $r(\overline{Q}^{r(s+1)}_{\alpha}) \le \frac{1}{2}2^{-l_{rs}}$. In both cases condition~3) is satisfied.

Suppose now that all cubes in this chain are blue.
Assume that
$\widetilde{Q}^{r(s+1)}_{\alpha} \bigcap \widetilde{Q}^{rs}_{\alpha'} \neq \varnothing$ and $r(Q^{r(s+1)}_{\alpha}) \geq \frac{1}{2} r(Q^{rs}_{\alpha'})$ for $\alpha' \in A^{rs}$, $\alpha \in A^{r(s+1)}$.
Then     $Q^{rs}_{\alpha'} \subset 8Q^{r(s+1)}_{\alpha}$, and hence,
$\widehat{\gamma}^{rs}_{\alpha'} \geq \frac{2}{q}\widehat{\gamma}^{r(s+1)}_{\alpha}$. On the other hand, by condition~1) of Definition~\ref{Def4.4} (as was pointed out above,
this condition is satisfied with $c_{1}=q^{3}$) and since all the cubes in the chain $\overline{Q}^{r(s+1)}_{\alpha} \subset \dots. \subset \overline{Q}^{rs}_{\alpha''}$ are blue and $r \geq 5$,
we have the estimate $\widehat{\gamma}^{r(s+1)}_{\alpha} \geq q^{4} \widehat{\gamma}^{rs}_{\alpha''}$ and hence $\widehat{\gamma}^{r(s+1)}_{\alpha} \geq q \widehat{\gamma}^{rs}_{\alpha'}$ (in view of condition 1)). This contradiction completes the verification of condition~3).

\textit{Step $2$}. We claim that estimate \eqref{eq4.1} holds. Arguing as in the proof of Lemma 3.1 of \cite{Tyu}, we see that
\begin{multline}
\label{eq4.4}
\int\limits_{\widetilde{Q}^{s}_{\alpha}}|\varphi^{s}_{\alpha}-\varphi(x)|\,dx \le\\
\le \frac{1}{|\widetilde{Q}^{s}_{\alpha}|}\int\limits_{\widetilde{Q}^{s}_{\alpha}}\int\limits_{\widetilde{Q}^{s}_{\alpha}}|\varphi(x)-\varphi(y)|\,dxdy \le
\iint\limits_{\widetilde{\Pi}^{s}_{\alpha}}|\nabla f(x,t)|\,dt dx, \qquad s \in \mathbb{Z}_{+}, \ \ \alpha \in \widetilde{A}^{s}.
\end{multline}

From \eqref{eq4.4} we see at once that
\begin{equation}
\label{eq4.5}
\begin{gathered}
\sum\limits_{s=0}^{\infty}\sum\limits_{\alpha \in \widetilde{A}^{s}} \widehat{\gamma}^{s}_{\alpha}\int\limits_{\widetilde{Q}^{s}_{\alpha}}|\varphi^{s}_{\alpha}-\varphi(x)|\,dx \le \sum\limits_{s=0}^{\infty}\sum\limits_{\alpha \in \widetilde{A}^{s}} \widehat{\gamma}^{s}_{\alpha}\iint\limits_{\widetilde{\Pi}^{s}_{\alpha}}|\nabla f(x,t)|\,dt\,tx\le \\
\le \sum\limits_{s=0}^{\infty}\sum\limits_{\substack{\alpha \in \widetilde{A}^{s}\\ \overline{Q}^{s}_{\alpha}\ \text{is yellow} }}\hbox{  }\widehat{\gamma}^{s}_{\alpha}\iint\limits_{\widetilde{\Pi}^{s}_{\alpha}}|\nabla f(x,t)|\,dt dx+\sum\limits_{s=0}^{\infty}\sum\limits_{\substack{\alpha \in \widetilde{A}^{s}\\ \overline{Q}^{s}_{\alpha}
\ \text{is blue} }} \hbox{  } \widehat{\gamma}^{s}_{\alpha} \iint\limits_{\widetilde{\Pi}^{s}_{\alpha}}|\nabla f(x,t)|\,dt\,tx= S_{1}+S_{2}.
\end{gathered}
\end{equation}

The sum $S_{1}$ is easily estimated by \eqref{eq4.2}. Using the finite (independent of  $j$ and~$m$) overlapping multiplicity of the sets $\widetilde{\Pi}_{l_{j},m}$ (when index $j$ is fixed and $m$ is variable) and Lemma~\ref{Lm2.1},
we arrive at the estimate
\begin{equation}
\begin{split}
\label{eq4.6}
&S_{1} \le C\sum\limits_{s=0}^{\infty}\sum\limits_{\substack{\alpha \in \widetilde{A}^{s}\\ \overline{Q}^{s}_{\alpha} \ \text{is yellow} }} \essinf \limits_{(x,t) \in \widetilde{\Pi}^{s}_{\alpha}}\gamma(x,t) \iint\limits_{\widetilde{\Pi}^{s}_{\alpha}}|\nabla f(x,t)|\,dt\,tx\le\\
&\le C \sum\limits_{j=0}^{\infty}\sum\limits_{m \in \mathbb{Z}^{n}} \iint\limits_{\widetilde{\Pi}_{l_{j},m}}\gamma(x,t)|\nabla f(x,t)|\,dt\,tx\le C \sum\limits_{j=0}^{\infty}\|f|W^{1}_{1}(\mathbb{R}^{n} \times (0,2^{-l_{j}}),\gamma)\| \le \\
&\le C \|f|W^{1}_{1}(\mathbb{R}^{n+1}_{+},\gamma)\|.
\end{split}
\end{equation}

We note that the constant $C > 0$ on the right of \eqref{eq4.6} depends only on $\lambda,n,C_{\gamma}$.

Given $s \in \mathbb{Z}_{+}$, $\alpha \in \widetilde{A}^{s}$, we set
$G^{s}_{\alpha}:=\widetilde{\Pi}^{s}_{\alpha} \setminus (\bigcup\limits_{\alpha' \in \widetilde{A}^{s+1}}\widetilde{\Pi}^{s}_{\alpha'}).$

The sum $S_{2}$ is estimated from above as follows (we change the order of summation in $(s,\alpha)$ and $(j,\alpha')$, respectively)
\begin{equation}
\begin{split}
\label{eq4.7}
&S_{2}\le \sum\limits_{s=0}^{\infty}\sum\limits_{\substack{\alpha \in \widetilde{A}^{s}\\ \overline{Q}^{s}_{\alpha} \ \text{is blue}}}\widehat{\gamma}^{s}_{\alpha}\sum\limits_{j=s}^{\infty}\sum\limits_{\substack{\alpha' \in \widetilde{A}^{j}\\ \widetilde{Q}^{s}_{\alpha}\bigcap \widetilde{Q}^{j}_{\alpha'} \neq \emptyset}}\iint\limits_{G^{j}_{\alpha'}}|\nabla f(x,t)|\,dxdt\\
&\le \sum\limits_{j=0}^{\infty}\sum\limits_{\alpha' \in \widetilde{A}^{j}} \sum\limits_{s=0}^{j}\sum\limits_{\substack{\alpha \in \widetilde{A}^{s}\\ \widetilde{Q}^{s}_{\alpha}\bigcap \widetilde{Q}^{j}_{\alpha'} \neq \emptyset \\ \overline{Q}^{s}_{\alpha} \ \text{is blue}}} \widehat{\gamma}^{s}_{\alpha} \|f|W^{1}_{1}(G^{j}_{\alpha'}\bigcap \widetilde{\Pi}^{s}_{\alpha})\|\Bigr).
\end{split}
\end{equation}

The main idea to be used for continuation of estimate \eqref{eq4.7} is close to that of \eqref{eq3.13}. However, here we are facing some substantial technical challenges.
First, the diameters of the sets $\widetilde{\Pi}^{s}_{\alpha}$ (from the right of \eqref{eq4.7}) may greatly differ from each other for a~fixed~$s$ and variable~$\alpha$.
Hence, the number of cubes $\widetilde{Q}^{j}_{\alpha'} \bigcap \widetilde{Q}^{s}_{\alpha} \neq \emptyset$ may be fairly large.
Second, a~more refined analysis of the behaviour of numbers $\widehat{\gamma}^{s}_{\alpha}$ is required.
We are unable to work with such numbers as with elements of  a~geometric progression (as this was done in \eqref{eq3.13}). Indeed,
taking numbers $\widehat{\gamma}^{s_{1}}_{\alpha_{1}}$ ($\alpha_{1} \in \widetilde{A}^{s_{1}}$) and $\widehat{\gamma}^{s_{2}}_{\alpha_{2}}$ ($\alpha_{2} \in \widetilde{A}^{s_{2}}$)
so as to have $\widetilde{Q}^{s_{1}}_{\alpha_{1}}\bigcap \widetilde{Q}^{j}_{\alpha'} \neq \emptyset$ and $\widetilde{Q}^{s_{2}}_{\alpha_{2}}\bigcap \widetilde{Q}^{j}_{\alpha'} \neq \emptyset$
we may not guarantee that at least one of the embeddings $Q^{s_{1}}_{\alpha_{1}} \subset Q^{s_{2}}_{\alpha_{2}}$ or $Q^{s_{2}}_{\alpha_{2}} \subset Q^{s_{1}}_{\alpha_{1}}$ hold.
The main idea will be to build a~chain of cubes
$\overline{Q}^{j}_{\alpha'} \subset ... \subset \overline{Q}^{0}_{\alpha''}$ and take care only about the numbers $\widehat{\gamma}^{j}_{\alpha'},..,\widehat{\gamma}^{0}_{\alpha''}$.
These numbers will play the role of a~`skeleton' which `supports' the remaining numbers $\widehat{\gamma}^{s}_{\alpha}$. Next, we split the set of numbers $\{\widehat{\gamma}^{j}_{\alpha'},..,\widehat{\gamma}^{0}_{\alpha''}\}$ into two sets. The first one will contain the numbers which behave like a geometric progression.
The other set will contain the numbers which, broadly speaking, do not behave like a~geometric progression. In the second case estimate \eqref{eq4.2}
will again prove useful. The formal proof proceeds as follows.

To continue estimating \eqref{eq4.7} we need the following important observation. We fix indexes $j \in \mathbb{Z}_{+}$ and $\alpha' \in \widetilde{A}^{j}$.
Given $s \in \{0,\dotsc,j\}$, we let $\overline{Q}^{s}_{\beta_{s}(\alpha')}$ ($\beta_{s}(\alpha') \in A^{s}$) denote the unique dyadic cube from the tiling $T^{s}$ that contains the cube $Q^{j}_{\alpha'}$.

We next use the fact that the system of  tilings $T=\{T^{s}\}$ is admissible for the weight~$\gamma$ (assertion~1)), apply assertion~2) of Lemma~ \ref{Lm4.1}, and
finally employ Lemma~\ref{Lm2.1}.  (For each fixed~$s$ and variable $\alpha$, the overlapping multiplicity of the sets $\widetilde \Pi^s_{\alpha}$ is finite and independent of~$s$ and~$\alpha$.)
We have
\begin{gather}
\notag
  \sum\limits_{s=0}^{j} \sum\limits_{\substack{\alpha \in \widetilde{A}^{s}\\ \widetilde{Q}^{s}_{\alpha}\bigcap \widetilde{Q}^{j}_{\alpha'} \neq \varnothing\\ \overline{Q}^{s}_{\alpha} \hbox{ is blue }}} \widehat{\gamma}^{s}_{\alpha} \|f|W^{1}_{1}(G^{j}_{\alpha'}\bigcap \widetilde{\Pi}^{s}_{\alpha})\|\le  c_{1}
\sum\limits_{s=0}^{j} \widehat{\gamma}^{s}_{\beta_{s}(\alpha')} \sum\limits_{\substack{\alpha \in \widetilde{A}^{s}\\ \widetilde{Q}^{s}_{\alpha}\bigcap \widetilde{Q}^{j}_{\alpha'} \neq \varnothing}}  \|f|W^{1}_{1}(G^{j}_{\alpha'}\bigcap \widetilde{\Pi}^{s}_{\alpha})\| \le \\
\le C(c_{1},n) \sum\limits_{s=0}^{j} \widehat{\gamma}^{s}_{\beta_{s}(\alpha')} \|f|W^{1}_{1}(G^{j}_{\alpha'})\|.
\label{eq4.8}
\end{gather}

We next partition the index set $\{0,\dots,j\}$ into two disjoint sets: $\{0,\dots,j\}={^{1}\Gamma^{j}_{\alpha}}\bigcup ^{2}\Gamma^{j}_{\alpha}$, where
$$
{}^{1}\Gamma^{s}_{\alpha'}:=\{s=0,\dotsc,j\mid  \overline{Q}^{s}_{\beta_{s}(\alpha)} \text{is blue}\}, ^{2}\Gamma^{j}_{\alpha'}:=\{s=0,\dotsc,j \mid \overline{Q}^{s}_{\beta_{s}(\alpha)} \ \text{is yellow}\}.
$$

We continue with estimate \eqref{eq4.7}. Using \eqref{eq4.8}, we have
\begin{multline}
\label{eq4.9}
S_{2} \le \sum\limits_{j=0}^{\infty}\sum\limits_{\substack{\alpha' \in \widetilde{A}^{j}}}(\sum\limits_{s \in ^{1}\Gamma^{j}_{\alpha'}}\widehat{\gamma}^{s}_{\beta_{s}(\alpha')})\|f|W^{1}_{1}(G^{j}_{\alpha'})\|+
\\
+\sum\limits_{j=0}^{\infty}\sum\limits_{\substack{\alpha' \in \widetilde{A}^{j}}}(\sum\limits_{s \in ^{2}\Gamma^{j}_{\alpha'}}\widehat{\gamma}^{j}_{\beta_{j}(\alpha')})\|f|W^{1}_{1}(G^{j}_{\alpha'})\|=:S_{2,1}+S_{2,2}.
\end{multline}

The following estimate is clear from the construction of the blue cubes:
\begin{equation}
\label{eq4.10}
\sum\limits_{s \in ^{1}\Gamma^{j}_{\alpha'}}\widehat{\gamma}^{s}_{\beta_{s}(\alpha')} \le q \widehat{\gamma}^{j}_{\alpha'}.
\end{equation}

From \eqref{eq4.10}, using Lemma \ref{Lm2.1} (here we use the finite overlapping multiplicity of the sets $G^{j}_{\alpha'}$, which is  independent of~$j$ and~$\alpha'$)
and  \eqref{eq2.4},~\eqref{eq2.5}, we get
\begin{multline}
\label{eq4.11}
S_{2,1} \le C \sum\limits_{j=0}^{\infty}\sum\limits_{\substack{\alpha' \in \widetilde{A}^{j}}} \widehat{\gamma}^{j}_{\alpha'}\|f|W^{1}_{1}(G^{j}_{\alpha'})\| \le
\\ \le
C \sum\limits_{j=0}^{\infty}\sum\limits_{\substack{\alpha' \in \widetilde{A}^{j}}} \|f|W^{1}_{1}(G^{j}_{\alpha'}, \gamma)\| \le C \|f|W^{1}_{1}(\mathbb{R}^{n+1}_{+},\gamma)\|.
\end{multline}

Given $i \in \mathbb{Z}_{+}$, $m \in \mathbb{Z}^{n}$, we set $E_{l_{i},m}:=\Pi_{l_{i},m} \setminus \bigcup\limits_{m' \in \mathbb{Z}^{n}}\Pi_{l_{i+1},m'}$.
For further purposes it is useful to recall that a~yellow cube is of the form  $\overline{Q}_{l_{j},m}$ with some $j \in \mathbb{Z}_{+}$ and $m \in \mathbb{Z}^{n}$.

To estimate $S_{2,2}$ we shall require the following key observation. We fix indexes $j \in \mathbb{Z}_{+}$ and $\alpha' \in \widetilde{A}^{j}$.
Let $E_{l_{i},m} \bigcap G^{j}_{\alpha'} \neq \varnothing$ for some $i \in \mathbb{Z}_{+}$, $m \in \mathbb{Z}^{n}$.
By elementary geometric considerations we see that $r(Q_{l_{i},m}) \le r(Q^{s_{0}}_{\beta_{s_{0}}(\alpha')})$, where $s_{0}=\max\{s|s\in ^{2}\Gamma^{j}_{\alpha'}\}$.
Indeed, otherwise we would have  $r(Q_{l_{i},m}) > r(Q^{s_{0}}_{\beta_{s_{0}}(\alpha')})$, and hence,
$r(Q_{l_{i+1},m}) \geq r(Q^{s_{0}}_{\beta_{s_{0}}(\alpha')})$. But then $G^{j}_{\alpha'} \subset \bigcup\limits_{m' \in \mathbb{Z}^{n}}\Pi_{l_{i+1},m'}$, which
contradicts the condition $E_{l_{i},m} \bigcap G^{j}_{\alpha'} \neq \emptyset$.

Moreover, if the cube $Q^{j}_{\alpha'} \subset Q^{s}_{\beta_{s}(\alpha')}$ with $s \in  ^{2}\Gamma^{j}_{\alpha'}$, then the cube $Q_{l_{i},m} \subset Q'^{s}_{\beta_{s}(\alpha')}$.
Here, the dyadic cube $Q'^{s}_{\beta_{s}(\alpha')}$ has common boundary points with the cube $Q^{s}_{\beta_{s}(\alpha')}$, and besides
$r(Q^{s}_{\beta_{s}(\alpha')})=r(Q'^{s}_{\beta_{s}(\alpha')})$. But this in combination with \eqref{eq.q'} implies that
\begin{equation}
\label{eq4.12}
\sum\limits_{s \in ^{2}\Gamma^{j}_{\alpha'}}\widehat{\gamma}^{s}_{\beta_{s}(\alpha')} \le q \sum\limits_{s=0}^{i}\sum\limits_{\substack{m' \in \mathbb{Z}^{n} \\ Q_{l_{i},m} \subset Q_{l_{s},m'}}} \widehat{\gamma}_{l_{s},m'}=:C t_{l_{i},m}.
\end{equation}

Using \eqref{eq4.12} and Lemma \ref{Lm2.1} (here we use the finite overlapping multiplicity of the sets $G^s_\alpha$, which is  independent of~$s$ and~$\alpha$),
we two times change the order of summation (first, with respect to $(j,\alpha')$ and $(i,m)$, and then with respect to $(i,m)$ and $(s,m')$) and take into account the equality
$\Pi_{l_{j},m'}=\bigcup\limits_{\substack{(i,m)\\Q_{l_{i},m} \subset Q_{l_{j},m'}}}E_{l_{i},m}$, estimate \eqref{eq2.1} and estimate~\eqref{eq4.2}. As a~result, we have
$$
S_{2,2} \le C \sum\limits_{j=0}^{\infty}\sum\limits_{\substack{\alpha' \in \widetilde{A}^{'}}} \sum\limits_{i=0}^{\infty}\sum\limits_{m \in \mathbb{Z}^{n}}t_{l_{i},m}\|f|W^{1}_{1}(G^{j}_{\alpha'}\bigcap E_{l_{i},m})\| \le
$$

$$
\le C \sum\limits_{i=0}^{\infty}\sum\limits_{m \in \mathbb{Z}^{n}}\biggl(\sum\limits_{s=0}^{i}\sum\limits_{\substack{m' \in \mathbb{Z}^{n} \\ Q_{l_{i},m} \subset Q_{l_{s},m'}}}\widehat{\gamma}_{l_{s},m'}\biggr) \|f|W^{1}_{1}(E_{l_{i},m})\| \le
$$

\begin{equation}
\begin{split}
\label{eq4.13}
& \le C \sum\limits_{s=0}^{\infty}\sum\limits_{m' \in \mathbb{Z}^{n}}\widehat{\gamma}_{l_{s},m'}\|f|W^{1}_{1}(\Pi_{l_{s},m'})\| \le \\
&\le C \sum\limits_{s=0}^{\infty}\sum\limits_{m' \in \mathbb{Z}^{n}}\|f|W^{1}_{1}(\Pi_{l_{s},m'},\gamma)\| \le C \sum\limits_{s=0}^{\infty}\|f|W^{1}_{1}(\mathbb{R}^{n}\times (0,2^{-l_{s}}),\gamma)\| \le \\ &\le C \|f|W^{1}_{1}(\mathbb{R}^{n+1}_{+},\gamma)\| .
\end{split}
\end{equation}

Combining estimates \eqref{eq4.9}, \eqref{eq4.11}, \eqref{eq4.13},  we find that
\begin{equation}
\label{eq4.14}
\begin{split}
S_{2} \le C \|f|W^{1}_{1}(\mathbb{R}^{n+1}_{+},\gamma)\|,
\end{split}
\end{equation}
where the constant $C$ depends only on $C_{\gamma},n,\lambda,c_{1},c_{2},q$.

Now estimate \eqref{eq4.1} follows from \eqref{eq3.6}, \eqref{eq4.5},  \eqref{eq4.6}, \eqref{eq4.14}.
This completes the proof of the theorem.

For further purposes we shall require a special partition of unity on $\mathbb{R}^{n} \times (0,2)$. Let $T=\{T^{s}\}^{\infty}_{s=0}(c_{1},c_{2})$
be a~system
of tilings of the space~$\mathbb{R}^{n}$ that is admissible for the weight~ $\gamma$.
The subsequent arguments will be carried out for $\lambda=2$,
even though they hold with minor technical modifications in the general case $\lambda=1+2^{-k}$, $k \in \mathbb{Z}_{+}$.
Hence, in what follows, $\widetilde{Q}_{k,m}=2Q_{k,m}$ for $(k,m) \in \mathbb{Z}_{+} \times \mathbb{Z}^{n}$ and
$\widetilde{Q}^{s}_{\alpha}=2Q^{s}_{\alpha}$ for $s \in \mathbb{Z}_{+}$, $\alpha \in \widetilde{A}^{s}$.

Given $k \in \mathbb{Z}_{+}$, $m \in \mathbb{Z}^{n}$, assume that  a~function
$\theta_{k,m} \in C^{\infty}_{0}(\mathbb{R}^{n})$ is such that  $\theta_{k,m}(x) \in (0,1]$ for $x  \in \widetilde{Q}_{k,m}$, $\theta_{k,m}(x)=0$ for
$x \in \mathbb{R}^{n} \setminus \widetilde{Q}_{k,m}$, and $|\nabla \theta_{k,m}(x)| \le C_{\theta} 2^{k}$ for $x \in \mathbb{R}^{n}$
with constant $C_{\theta} > 0$ independent both of~ $k$, ~$m$ and $T$. We also assume that  $\sum\limits_{m \in \mathbb{Z}^{n}}\theta_{k,m} \equiv 1$ on~$\mathbb{R}^{n}$.

Next, given $k \in \mathbb{Z}_{+}$,  assume that  a function $\psi_{k} \in C^{\infty}_{0}((0,\infty))$ is such that
$\psi_{k}(t) \in (0,1)$ for $t \in (\frac{7}{8}2^{-k},\frac{9}{8}2^{-k+1})$, $\psi_{k}(t)=0$ for
$t \in  (0,+\infty) \setminus (\frac{7}{8}2^{-k},\frac{9}{8}2^{-k+1})$, and $\biggl|\frac{d\psi_{k}}{dt}(t)\biggr| \le C_{\psi}2^{k}$ for $t > 0$ with constant $C_{\psi} > 0$ independent both of~ $s$ and~$\alpha$. We also assume that  $\sum\limits_{k \in \mathbb{Z}_{+}}\psi_{k} \equiv 1$ on~$(0,2)$.

  We set $\Theta_{k,m}=\theta_{k,m}\psi_{k}$ for $k \in \mathbb{Z}_{+}$, $m \in \mathbb{Z}^{n}$. It is clear that $\Theta_{k,m} \in C^{\infty}_{0}(\mathbb{R}^{n+1}_{+})$ and

\begin{equation}
\label{eq4.16}
\sum\limits_{k \in \mathbb{Z}_{+}}\sum\limits_{m \in \mathbb{Z}^{n}}\Theta_{k,m}(x,t)=1, \qquad (x,t) \in \mathbb{R}^{n} \times (0,2).
\end{equation}

For every $(k,m) \in \mathbb{Z}_{+} \times \mathbb{Z}^{n}$ we have

\begin{equation}
\label{eq4.17}
|\nabla\Theta_{k,m}(x,t)| \le C_{\theta}C_{\psi} 2^{k}, \qquad (x,t) \in \mathbb{R}^{n+1}_{+}.
\end{equation}

In what follows we shall require some combinatoric arguments. Recall, that we are dealing with the case $\lambda=2$, and so
$\widetilde{Q}^{s}_{\alpha}=2Q^{s}_{\alpha}$. Given $s \in \mathbb{Z}_{+}$, $\alpha \in \widetilde{A}^{s}$, we set $B^{s}_{\alpha}:=\{(k,m) \in \mathbb{Z}_{+} \times \mathbb{Z}^{n}: Q_{k,m} \subset \widetilde{Q}^{s}_{\alpha}\}$.

For any fixed $s \in \mathbb{Z}_{+}$, we represent the index set  $\widetilde{A}^{s}$ as a union of finite number (by condition~4) of Definition \ref{Def4.4})
of at most countable index subsets $\widetilde{A}^{s,k}$, $k \in \{1,...,t(s)\}$, which are pairwise disjoint. Besides, we shall require that
$r(Q^{s}_{\alpha})=r(Q^{s}_{\alpha'})$ for $\alpha,\alpha' \in \widetilde{A}^{s,k}$ ($k \in \{1,...,t(s)\}$) and $r(Q^{s}_{\alpha}) < r(Q^{s}_{\alpha'})$ for $\alpha \in \widetilde{A}^{s,k}$, $\alpha' \in \widetilde{A}^{s,k+1}$ ($k \in \{1,...,t(s)-1\}$).

Now, given fixed $s \in \mathbb{Z}_{+}$ and $k \in \{1,...,t(s)\}$, we label the cubes
$\{\widetilde{Q}^{s}_{\alpha}\}_{\alpha \in \widetilde{A}^{s,k}}$ by natural number; that is, $\{\widetilde{Q}^{s}_{\alpha}\}_{\alpha \in
\widetilde{A}^{s,k}}=\{\widetilde{Q}^{s}_{\alpha_{i}}\}_{i=1}^{\infty}$ (for each $s$ and~$k$ the procedure of labeling is, in general, different, but for us this is
immaterial). Let $D^{s}_{\alpha_{1}}=B^{s}_{\alpha_{1}}$. If for some $k' \in \mathbb{N}$ we have already constructed the index sets $D^{s}_{\alpha_{j}}$ ($j \in \{1,..,k'\}$),
then we set $D^{s}_{\alpha_{k'+1}}:=B^{s}_{\alpha_{k'+1}} \setminus \bigcup\limits_{j=1}^{k'}D^{s}_{\alpha_{j}}$.
So, by induction, for each $k' \in \mathbb{N}$ we construct the index set $D^{s}_{\alpha_{k'}}$.
Arguing similarly for all $k \in \{1,...,t(s)\}$ and next for all $s \in \mathbb{Z}_{+}$, we shall construct the index sets $D^{s}_{\alpha}$ for any $s \in \mathbb{Z}_{+}$ and any
$\alpha \in \widetilde{A}^{s}$. Finally, we set
$$
E^{s}_{\alpha}:=D^{s}_{\alpha} \setminus \bigcup\limits_{\substack{s' \geq s, \alpha' \in \widetilde{A}^{s'} \\ r(Q^{s'}_{\alpha'}) < r(Q^{s}_{\alpha})}}B^{s'}_{\alpha'}.
$$

Note that the sets $D^{s}_{\alpha}$ for $s \in \mathbb{Z}_{+}$, $k \in \{1,..,t(s)\}$, $\alpha \in \widetilde{A}^{s,k}$ are pairwise disjoint.
Hence it clearly follows from the inclusion $D^{s}_{\alpha} \subset B^{s}_{\alpha}$ that $E^{s}_{\alpha} \bigcap E^{s'}_{\alpha'} = \emptyset$ for $(s,\alpha) \neq (s',\alpha')$.
It is easily checked that $B^{s}_{\alpha} \subset \bigcup\limits_{\substack{s' \geq s, \alpha' \in \widetilde{A}^{s'} \\ r(Q^{s'}_{\alpha'}) \le r(Q^{s}_{\alpha})}}D^{s'}_{\alpha'}$ for any
$s \in \mathbb{Z}_{+}$, $\alpha \in \widetilde{A}^{s}$. Hence, from the definition of the sets $E^{s}_{\alpha}$ and conditions~3) of Definition \ref{Def4.4} one readily verifies that
\begin{equation}
\label{eq4.18}
\bigcup\limits_{s \in \mathbb{Z}_{+}}\bigcup\limits_{\alpha \in \widetilde{A}^{s}}E^{s}_{\alpha}=\mathbb{Z}_{+}\times \mathbb{Z}^{n}.
\end{equation}

We set
\begin{equation}
\label{eq4.19}
g^{s}_{\alpha}(x,t):=\sum\limits_{(k,m) \in E^{s}_{\alpha}}\Theta_{k,m}(x,t), \qquad (x,t) \in \mathbb{R}^{n+1}_{+}.
\end{equation}

The next lemma follows from \eqref{eq4.18} and \eqref{eq4.19}.

\begin{Lm}
\label{Lm4.2}
The functions $g^{s}_{\alpha}$ have the following properties:
\begin{list}{}{\itemsep=0pt\topsep=2pt\parsep=0pt}
\item[\rm  1)] $g^{s}_{\alpha} \in C^{\infty}_{0}(\mathbb{R}^{n+1}_{+})$ for $s \in \mathbb{Z}_{+}$, $\alpha \in \widetilde{A}^{s}$,
\item[\rm  2)] $\sum\limits_{s=0}^{\infty}\sum\limits_{\alpha \in \widetilde{A}^{s}}g^{s}_{\alpha}(x,t)=1$ for $(x,t) \in \mathbb{R}^{n} \times (0,2)$,
\item[\rm  3)] for any point $(x,t) \in \mathbb{R}^{n+1}_{+}$ there exist at most $C(n)$ functions $g^{s}_{\alpha}$ for which
$g^{s}_{\alpha}(x,t) > 0$,
\item[\rm  4)] for any $s \in \mathbb{Z}_{+}$, $\alpha \in \widetilde{A}^{s}$
\begin{equation}
\label{eq4.20}
|\nabla g^{s}_{\alpha}(x,t)| \le C 2^{k}, \qquad (x,t) \in \overline{Q}_{k,m} \times [\frac{1}{2^{k}},\frac{1}{2^{k-1}}].
\end{equation}
\end{list}
The constant $C > 0$ on the right of \eqref{eq4.20} depends only on $n, C_{\psi}, C_{\theta}$.
\end{Lm}

 Henceforward, $\mu_{n}$ will denote the Lebesgue measure in $\mathbb{R}^{n}$.

The following fact will be crucial to all our subsequent work.

\begin{Lm}
\label{Lm4.3}
Let $\gamma \in A^{loc}_{1}(\mathbb{R}^{n+1})$, $\lambda=2$, $c_{1},c_{2} > 0$. Let $T=\{T^{s}\}^{\infty}_{s=0}(c_{1},c_{2})$
be a~system
of tilings of the space~$\mathbb{R}^{n}$ that is admissible for the weight~ $\gamma$. Then for every $s \in \mathbb{Z}_{+}$, $\alpha \in \widetilde{A}^{s}$ the following inequality holds
\begin{equation}
\label{eq4.21}
\iint\limits_{\operatorname{supp}g^{s}_{\alpha}}\gamma(x,t)|\nabla g^{s}_{\alpha}(x,t)|\,dxdt \le C \widehat{\gamma}^{s}_{\alpha}\mu_{n}(Q^{s}_{\alpha}).
\end{equation}

The constant $C > 0$ depends only on $n, c_{1},c_{2}, C_{\psi}, C_{\theta}, C_{\gamma}$.
\end{Lm}

\textbf{Proof.} Given any $(k,m) \in \mathbb{Z}_{+} \times \mathbb{Z}^{n}$, let  $\hat{Q}_{k,m}:=Q_{k,m} \times (\frac{1}{2^{k}},\frac{1}{2^{k-1}})$.
For fixed $s \in \mathbb{Z}_{+}$ and $\alpha \in \widetilde{A}^{s}$ we consider only those cubes  $\hat{Q}_{k,m}$, $(k,m) \in E^{s}_{\alpha}$, for which the function
$g^{s}_{\alpha}$ is not identically zero on  $2\hat{Q}_{k,m}$. Clearly, the number of such cubes is finite (in general, depending on $s$ and~$\alpha$). Let
$\{Q_{j}\}_{j=1}^{n(s,\alpha)}$ be the above set of cubes. By \eqref{eq4.19}, \eqref{eq4.20} we have, for any $s \in \mathbb{Z}_{+}$ and $\alpha \in \mathbb{Z}^{n}$,
$$
|\nabla g^{s}_{\alpha}(x,t)| \le C(n,C_{\psi},C_{\theta})(r(Q_{j}(s,\alpha)))^{-1}, \qquad (x,t) \in 2Q_{j}(s,\alpha).
$$

Hence,
\begin{equation}
\label{eq4.22}
\iint\limits_{\operatorname{supp}g^{s}_{\alpha}}\gamma(x,t)|\nabla g^{s}_{\alpha}(x,t)|\,dxdt \le C\sum\limits_{j=1}^{n(s,\alpha)}\mu_{n}(\check{Q}_{j}(s,\alpha))\frac{1}{\mu_{n+1}(Q_{j}(s,\alpha))}\iint\limits_{2Q_{j}(s,\alpha)}\gamma(x,t)\,dxdt.
\end{equation}

We shall henceforward denote by $\check{Q}_{j}(s,\alpha)$ the projections of the cube  $Q_{j}(s,\alpha)$ to the hyperplane $\mathbb{R}^{n} \times \{0\}$.

Note that, for any $j \in \{1,..,n(s,\alpha)\}$, the side length $r(Q_{j}(s,\alpha)) \geq r(Q^{s+1}_{\alpha'})$ for some $\alpha' \in \widetilde{A}^{s+1}$ for which
$\check{Q}_{j}(s,\alpha) \bigcap \widetilde{Q}^{s+1}_{\alpha'} \neq \emptyset$. Indeed, otherwise $r(Q_{j}(s,\alpha)) \le 2r(Q^{s+1}_{\alpha'})$ for all $\alpha' \in \widetilde{A}^{s+1}$
for which $\check{Q}_{j}(s,\alpha) \bigcap \widetilde{Q}^{s+1}_{\alpha'} \neq \emptyset$. Hence, ${Q}_{j}(s,\alpha) \subset \bigcup\limits_{\alpha' \in \widetilde{A}^{s+1}}
\bigcup\limits_{(k,m) \in B^{s+1}_{\alpha'}}\overline{\hat {Q}}_{k,m}$, which shows that the cube $Q_{j}(s,\alpha)$ cannot be contained in the set $g^{s}_{\alpha}$.
But this contradicts the construction of the cubes $\{Q_{j}(s,\alpha)\}_{j=1}^{n(s,\alpha)}$.

Thus, from the above we have $Q^{s+1}_{\alpha'} \times (0,r(Q^{s+1}_{\alpha'})) \subset 8Q_{j}(s,\alpha) \subset 8\widetilde{Q}^{s}_{\alpha} \times (0,r(Q^{s}_{\alpha}))$.
Hence, using \eqref{eq4.1}, \eqref{eq4.2} and conditions 1), 2) of Definition~\ref{Def4.4},
\begin{equation}
\label{eq4.23}
\frac{1}{C(q,n,c_{1},c_{2})} \widehat{\gamma}^{s}_{\alpha} \le  \frac{1}{\mu_{n+1}(Q_{j}(s,\alpha))}\iint\limits_{2Q_{j}(s,\alpha)}\gamma(x,t)\,dxdt \le C(q,n,c_{1},c_{2})\widehat{\gamma}^{s}_{\alpha}.
\end{equation}

From \eqref{eq4.22}, \eqref{eq4.23} we conclude that the lemma will be proved once we prove the estimate
\begin{equation}
\label{eq4.24}
\sum\limits_{j=1}^{n(s,\alpha)}\mu_{n}(\check{Q}_{j}(s,\alpha)) \le C \mu_{n}(Q^{s}_{\alpha}),
\end{equation}
in which the constant $C > 0$ depends only on $n$.

We fix indexes $s \in \mathbb{Z}_{+}$, $\alpha \in \widetilde{A}^{s}$ and a number $l \in \mathbb{N}$. Consider the set
$$
U(s,\alpha,l):=\widetilde{Q}^{s}_{\alpha} \setminus \bigcup\limits_{s'=s}^{s+1}\bigcup\limits_{\substack{\alpha' \in \widetilde{A}^{s'} \\ \widetilde{Q}^{s}_{\alpha}\bigcap \widetilde{Q}^{s'}_{\alpha'} \neq \emptyset \\ 2^{-l}r(Q^{s}_{\alpha})\le r(Q^{s'}_{\alpha'}) < r(Q^{s}_{\alpha})}}\widetilde{Q}^{s'}_{\alpha'}.
$$

It is easily checked that $U(s,\alpha,l+1) \subset U(s,\alpha,l) \subset \widetilde{Q}^{s}_{\alpha}$ for $l \in \mathbb{N}$.

Note that if the function $\nabla g^{s}_{\alpha}$ is not identically zero on the cube $2Q_{j}(s,\alpha)$ with side length $r(Q_{j}(s,\alpha))=2^{-l}r(Q^{s}_{\alpha})$ for $l \in \mathbb{N}$, then
\begin{equation}
\label{eq4.41'}
\check{Q}_{j}(s,\alpha) \bigcap \partial(U(s,\alpha,l)\bigcup U(s,\alpha,l+1)) \neq \emptyset.
\end{equation}

It is also easy to see that for  $l \in \mathbb{N}$
\begin{equation}
\label{eq4.25}
\sum\limits_{\substack{j\in \{1,..,n(s,\alpha)\} \\ r(Q_{j}(s,\alpha))=2^{-l}r(Q^{s}_{\alpha}) \\ 2Q_{j}(s,\alpha) \bigcap \partial U(s,\alpha,l+1)\neq \emptyset}}\mu_{n}(\check{Q}_{j}(s,\alpha)) \le C(n)\sum\limits_{\substack{j\in \{1,..,n(s,\alpha)\} \\ r(Q_{j}(s,\alpha))=2^{-l}r(Q^{s}_{\alpha}) \\ 2Q_{j}(s,\alpha) \bigcap \partial U(s,\alpha,l)\neq \emptyset}}\mu_{n}(\check{Q}_{j}(s,\alpha)) .
\end{equation}

 Next, we may assume that $n \geq 2$, for otherwise the arguments in the case $n=1$ are substantially easier.

The key observation is that, for every $l \in \mathbb{N}$,
\begin{equation}
\label{eq4.26}
\begin{split}
&\sum\limits_{\substack{j\in \{1,..,n(s,\alpha)\} \\ r(Q_{j}(s,\alpha))=2^{-l}r(Q^{s}_{\alpha}) \\ \check{Q}_{j}(s,\alpha) \bigcap \partial U(s,\alpha,l)\neq \emptyset}}\mu_{n}(\check{Q}_{j}(s,\alpha)) \le \sum\limits_{s'=s}^{s+1}\sum\limits_{\substack{\alpha' \in \widetilde{A}^{s} \\ \widetilde{Q}^{s'}_{\alpha'}\bigcap \widetilde{Q}^{s}_{\alpha} \neq \emptyset \\ 2^{-l}r(Q^{s}_{\alpha})\le r(Q^{s'}_{\alpha'}) < r(Q^{s}_{\alpha})}}2^{-l}r(Q^{s}_{\alpha})\mu_{n-1}(\partial(\widetilde{Q}^{s'}_{\alpha'}\bigcap \widetilde{Q}^{s}_{\alpha})) \le \\
&\le C(n) \sum\limits_{s'=s}^{s+1}\sum\limits_{\substack{\alpha' \in \widetilde{A}^{s} \\ \widetilde{Q}^{s'}_{\alpha'}\bigcap \widetilde{Q}^{s}_{\alpha} \neq \emptyset \\ 2^{-l}r(Q^{s}_{\alpha})\le r(Q^{s'}_{\alpha'}) < r(Q^{s}_{\alpha})}}\frac{2^{-l}r(Q^{s}_{\alpha})}{r(Q^{s'}_{\alpha'})}|\widetilde{Q}^{s'}_{\alpha'}\bigcap \widetilde{Q}^{s}_{\alpha}|.
\end{split}
\end{equation}

From \eqref{eq4.41'}, \eqref{eq4.25}, \eqref{eq4.26},
and taking into account that $\frac{2^{-l}r(Q^{s}_{\alpha})}{r(Q^{s'}_{\alpha'})}=2^{-j}$ (on the right of \eqref{eq4.26}), we have, for some $j \in \mathbb{Z}_{+}$,
\begin{equation}
\label{eq4.27}
\begin{split}
&\sum\limits_{j\in \{1,..,n(s,\alpha)\}} \mu_{n}(\check{Q}_{j}(s,\alpha)) \le C(n)\mu_{n}(Q^{s}_{\alpha}) +C(n) \sum\limits_{l=1}^{\infty}\sum\limits_{s'=s}^{s+1}\sum\limits_{\substack{\alpha' \in \widetilde{A}^{s} \\ \widetilde{Q}^{s'}_{\alpha'}\bigcap \widetilde{Q}^{s}_{\alpha} \neq \emptyset \\ 2^{-l}r(Q^{s}_{\alpha})\le r(Q^{s'}_{\alpha'}) < r(Q^{s}_{\alpha})}}\frac{2^{-l}r(Q^{s}_{\alpha})}{r(Q^{s'}_{\alpha'})}|\widetilde{Q}^{s'}_{\alpha'}\bigcap \widetilde{Q}^{s}_{\alpha}|\\
&\le C(n) \sum\limits_{s'=s}^{s+1}\sum\limits_{\substack{\alpha' \in \widetilde{A}^{s} \\ \widetilde{Q}^{s'}_{\alpha'}\bigcap \widetilde{Q}^{s}_{\alpha} \neq \emptyset}}\sum\limits_{j=1}^{\infty}2^{-j}|\widetilde{Q}^{s'}_{\alpha'}\bigcap \widetilde{Q}^{s}_{\alpha}| \le C(n)\mu_{n}(Q^{s}_{\alpha}).
\end{split}
\end{equation}

Now estimate \eqref{eq4.24} follows from \eqref{eq4.27}. The proof of the lemma is complete.

\begin{Th}
\label{Th4.2}
Let a weight $\gamma \in A^{\rm loc}_{1}(\mathbb{R}^{n+1})$, $c_{1},c_{2} \geq 1$. Assume that for a~function $\varphi \in L^{loc}_{1}(\mathbb{R}^{n})$ there
 exists a~system of tilings $T=\{T^{s}\}(c_{1},c_{2})$ admissible for the weight~$\gamma$  such that
$$
\sum\limits_{m \in \mathbb{Z}^{n}}\widehat{\gamma}_{0,m}\varphi_{0,m}+\sum\limits_{s=1}^{\infty}\sum\limits_{\alpha \in \widetilde{A}^{s}} \widehat{\gamma}^{s}_{\alpha}\int\limits_{\widetilde{Q}^{s}_{\alpha}}|\varphi^{s}_{\alpha}-\varphi(x)|\,dx < \infty
$$

Then there exists a~function $f \in W^{1}_{1}(\mathbb{R}^{n+1}_{+},\gamma)$ such that  $\varphi = \operatorname{tr}\left|_{t=0}f\right.$, and moreover,
\begin{equation}
\label{eq4.28}
\begin{split}
&C \|f|W^{1}_{1}(\mathbb{R}^{n+1}_{+},\gamma)\|
\le \sum\limits_{m \in \mathbb{Z}^{n}}\widehat{\gamma}_{0,m}\|\varphi|L_{1}(Q_{0,m})\|+ \sum\limits_{s=1}^{\infty}\sum\limits_{\alpha \in \widetilde{A}^{s}} \widehat{\gamma}^{s}_{\alpha}\int\limits_{\widetilde{Q}^{s}_{\alpha}}|\varphi^{s}_{\alpha}-\varphi(x)|\,dx .
\end{split}
\end{equation}
The constant $C > 0$ on the left of \eqref{eq4.28} depends only on $n,C_{\gamma},\lambda, c_{1},c_{2}$.
\end{Th}

 \textbf{Proof.} We shall prove the theorem for $\lambda=2$,
but our arguments will hold in the general case $\lambda = 1+2^{-k}$, $k \in \mathbb{Z}_{+}$ with minor technical modifications.

\textit{Step $1$}. We set
\begin{equation}
\label{eq4.29}
f(x,t)=\sum\limits_{s=0}^{\infty}\sum\limits_{\alpha \in \widetilde{A}^{s}} g^{s}_{\alpha}(x,t)\varphi^{s}_{\alpha}, \qquad (x,t) \in \mathbb{R}^{n+1}_{+}.
\end{equation}

Note that the function $f \in C^{\infty}(\mathbb{R}^{n+1}_{+})$. We claim that \eqref{eq4.28} holds.

To this aim we first estimate the integral
\begin{gather*}
J:=\iint\limits_{\mathbb{R}^{n+1}_{+}}\gamma(x,t)|\nabla f(x,t)|\,dx\,dt=\\
=\int\limits_{\mathbb{R}^{n}}\int\limits_{0}^{2}\gamma(x,t)|\nabla f(x,t)|\,dx\,dt+\int\limits_{\mathbb{R}^{n}}\int\limits_{2}^{\infty}\gamma(x,t)|\nabla f(x,t)|\,dx\,dt=:J_{1}+J_{2}.
\end{gather*}

From Lemma 4.2 we have
$$
  |\nabla f(x,t)|\le C \sum_{\alpha\in \tilde A^0} |\varphi ^0_\alpha| \chi_{\widetilde{Q}^0_\alpha}(x), \qquad x\in \mathbb R^n,\ \ t\ge 2.
$$

The cubes $\widetilde{Q}_{0,m}$ have finite (independent of~$m$) overlapping multiplicity, and hence by \eqref{eq2.3} we have
\begin{equation}
\label{eq4.30}
J_{2} \le C \sum\limits_{m \in \mathbb{Z}^{n}}\widehat{\gamma}_{0,m}\biggl(\sum\limits_{\substack{ m' \in \mathbb{Z}^{n} \\ \widetilde{Q}_{0,m'}\bigcap \widetilde{Q}_{0,m} \neq \varnothing}}|\varphi_{0,m'}| \biggr) \le C \sum\limits_{m \in \mathbb{Z}^{n}}\widehat{\gamma}_{0,m}\|\varphi|L_{1}(Q_{0,m})\|.
\end{equation}

Clearly, the constant $C$ in~\eqref{eq4.30} depends only on~$n,C_{\gamma}$.

Now let us estimate the more involved integral  $J_{1}$. Since $\mathbb{R}^{n} \times (0,2) \subset \bigcup\limits_{s \in \mathbb{Z}_{+}}\bigcup\limits_{\alpha \in \widetilde{A}^{s}}\operatorname{supp}g^{s}_{\alpha}$
 we find that
\begin{equation}
\label{eq4.31}
J_{1} \le  \sum\limits_{s=0}^{\infty}\sum\limits_{\alpha \in \widetilde{A}^{s}}\hbox{ }\iint\limits_{\operatorname{supp}g^{s}_{\alpha} \bigcap \mathbb{R}^{n} \times (0,2)}\gamma(x,t)|\nabla f(x,t)|\,dx\,dt.
\end{equation}

Given a fixed index $s_{0} \in \mathbb{Z}_{+}$ and $\alpha_{0} \in \widetilde{A}^{s_{0}}$, we use Lemma \ref{Lm4.2} (assertions~1), 2), 3)) and recall that
the system of tilings~$T$ is admissible (condition~3) of Definition~\ref{Def4.4}). We have (if $s_{0}=0$ we set formally $s_{0}-1=0$)
\begin{equation}
\begin{split}
\label{eq4.32}
&\iint\limits_{\operatorname{supp}g^{s_{0}}_{\alpha_{0}} \bigcap \mathbb{R}^{n}\times (0,2)}\gamma(x,t)|\nabla f(x,t)|\,dx\,dt = \iint\limits_{\operatorname{supp}g^{s_{0}}_{\alpha_{0}} \bigcap \mathbb{R}^{n}\times (0,2)}\gamma(x,t)\left|\sum\limits_{s=0}^{\infty}\sum\limits_{\alpha \in \widetilde{A}^{s}}\nabla g^{s}_{\alpha}(x,t)\varphi^{s}_{\alpha}\right|\,dx\,dt =\\
&= \iint\limits_{\operatorname{supp}g^{s_{0}}_{\alpha_{0}} \bigcap \mathbb{R}^{n}\times (0,2)}\gamma(x,t)\left|\sum\limits_{s=0}^{\infty}\sum\limits_{\alpha \in \widetilde{A}^{s}}\nabla g^{s}_{\alpha}(x,y)(\varphi^{s}_{\alpha}-\varphi^{s_{0}}_{\alpha_{0}})\right|\,dx\,dt \le \\
&\le \sum\limits_{s=s_{0}-1}^{s_{0}+1}\sum\limits_{\substack{\alpha \in \widetilde{A}^{s} \\ \operatorname{supp}g^{s_{0}}_{\alpha_{0}}\bigcap\operatorname{supp}g^{s}_{\alpha} \neq \varnothing}}   \left(\hbox{ }\iint\limits_{\operatorname{supp}g^{s_{0}}_{\alpha_{0}} \bigcap \mathbb{R}^{n}\times (0,2)}|\nabla g^{s}_{\alpha}(x,t)|\gamma(x,t)\,dx\,dt\right)|\varphi^{s}_{\alpha}-\varphi^{s_{0}}_{\alpha_{0}}|.
\end{split}
\end{equation}

The main crux now is to estimate  $|\nabla g^{s}_{\alpha}(x,t)|$ on the set $\operatorname{supp}g^{s_{0}}_{\alpha_{0}} \bigcap \mathbb{R}^{n}\times (0,2)$.
By Lemma \ref{Lm4.3} and using conditions 1), 2) of Definition~\ref{Def4.4}, we conclude that, for $s \in \{s_{0}-1,s_{0},s_{0}+1\}$,
\begin{equation}
\begin{split}
\label{eq4.33}
&\iint\limits_{\operatorname{supp}g^{s_{0}}_{\alpha_{0}} \bigcap \operatorname{supp}g^{s}_{\alpha}}\gamma(x,t)|\nabla g(x,t)|\,dx\,dt \le C \widehat{\gamma}^{s_{0}}_{\alpha_{0}} \min\{\mu_{n}(Q^{s}_{\alpha}), \mu_{n}(Q^{s_{0}}_{\alpha_{0}})\} \le \\
&\le C \widehat{\gamma}^{s_{0}}_{\alpha_{0}} \mu_{n}(\widetilde{Q}^{s}_{\alpha} \bigcap \widetilde{Q}^{s_{0}}_{\alpha_{0}}).
\end{split}
\end{equation}

Substituting estimate \eqref{eq4.33} into \eqref{eq4.32} and using conditions 1), 2) of Definition~ \ref{Def4.4}, this gives
\begin{equation}
\begin{split}
\label{eq4.34}
&\iint\limits_{\operatorname{supp}g^{s_{0}}_{\alpha_{0}} \bigcap \mathbb{R}^{n}\times (0,2)}\gamma(x,t)|\nabla f(x,t)|\,dx\,dt \le C \sum\limits_{s=s_{0}-1}^{s_{0}+1}\sum\limits_{\substack{\alpha \in \widetilde{A}^{s} \\ \operatorname{supp}g^{s_{0}}_{\alpha_{0}}\bigcap\operatorname{supp}g^{s}_{\alpha} \neq \varnothing}}   \hbox{ }\widehat{\gamma}^{s_{0}}_{\alpha_{0}} \mu_{n}(\widetilde{Q}^{s}_{\alpha} \bigcap \widetilde{Q}^{s_{0}}_{\alpha_{0}})|\varphi^{s}_{\alpha}-\varphi^{s_{0}}_{\alpha_{0}}| \le \\
&\le C \sum\limits_{s=s_{0}-1}^{s_{0}+1}\sum\limits_{\substack{\alpha \in \widetilde{A}^{s} \\ \operatorname{supp}g^{s_{0}}_{\alpha_{0}}\bigcap\operatorname{supp}g^{s}_{\alpha} \neq \varnothing}} \left( \widehat{\gamma}^{s}_{\alpha} \int\limits_{\widetilde{Q}^{s}_{\alpha}\bigcap \widetilde{Q}^{s_{0}}_{\alpha_{0}}}|\varphi(x)-\varphi^{s}_{\alpha}|\,dx +  \widehat{\gamma}^{s_{0}}_{\alpha_{0}} \int\limits_{\widetilde{Q}^{s}_{\alpha}\bigcap \widetilde{Q}^{s_{0}}_{\alpha_{0}}}|\varphi(x)-\varphi^{s_{0}}_{\alpha_{0}}|\,dx\right).
\end{split}
\end{equation}

Summing estimate  \eqref{eq4.34} over all indexes $s_{0}$, $\alpha_{0}$, taking into account conditions 1), 2) of Definition \ref{Def4.4}, using assertion~2) of Lemma~\ref{Lm4.1}
and employing Lemma~\ref{Lm2.1} (with $d=n$), we finally have
\begin{equation}
\begin{split}
\label{eq4.35}
J_{1} \le C(c_{1},c_{2},C_{1},C_{2},n) \biggl(\sum\limits_{s=1}^{\infty}\sum\limits_{\alpha \in \widetilde{A}^{s}} \widehat{\gamma}^{s}_{\alpha}\int\limits_{\widetilde{Q}^{s}_{\alpha}}|\varphi^{s}_{\alpha}-\varphi(x)|\,dx + \sum\limits_{m \in \mathbb{Z}^{n}} \widehat{\gamma}_{0,m}\|\varphi|L_{1}(Q_{0,m})\|\biggr).
\end{split}
\end{equation}

Arguing as in the estimate 3.11 of \cite{Tyu} we have

\begin{equation}
\label{eq4.36}
\iint\limits_{\mathbb{R}^{n+1}_{+}}\gamma(x,t)|f(x,t)|\,dxdt \le C \iint\limits_{\mathbb{R}^{n+1}_{+}}\gamma(x,t)|\nabla f(x,t)|\,dxdt.
\end{equation}

Now \eqref{eq4.28} follows from \eqref{eq4.30}, \eqref{eq4.35}, \eqref{eq4.36}.

\textit{Step $2$}. We now claim that $\varphi =\hbox{tr}|_{t=0}f$.

For any fixed $t \in (0,1)$ from assertions 2) of Lemma~\ref{Lm4.2}, from condition~3) of Definition \ref{Def4.4}, and from \eqref{eq4.19} we have the following estimate
\begin{equation}
\label{eq4.37}
|f(x,t)-\varphi(x)| = \left|\sum\limits_{s=s(t)-1}^{s=s(t)+1}\sum\limits_{\substack{\alpha \in \widetilde{A}^{s}\\ x \in \widetilde{Q}^{s}_{\alpha}}}g^{s}_{\alpha}(x,t)(\varphi^{s}_{\alpha}-\varphi(x))\right| \le \sum\limits_{s=s(t)-1}^{s=s(t)+1}\sum\limits_{\substack{\alpha \in \widetilde{A}^{s}\\ x \in \widetilde{Q}^{s}_{\alpha}}}\frac{1}{|\widetilde{Q}^{s}_{\alpha}|}\int\limits_{\widetilde{Q}^{s}_{\alpha}}|\varphi(\widetilde{x})-\varphi(x)|\,d\widetilde{x}.
\end{equation}

Note that the set $\tilde Q^s_\alpha$ of cubes containing the point $x$ forms a~regular family in the sense of \S\,1.8 of \cite{St2}.
Combining the arguments of \S\,1.8 of \cite{St2} with condition~3) from Definition~\ref{Def4.4} and taking into account the finite (depending only on~$n$) overlapping multiplicity of the cubes $\widetilde{Q}^{s}_{\alpha}$ (when $s$ is fixed and $\alpha$ variable) it is easily deduce from \eqref{eq4.37} that
\begin{equation}
\label{eq4.38}
\varphi(x)=\lim\limits_{t \to +0}f(x,t) \qquad \hbox{for almost all } x \in \mathbb{R}^{n}.
\end{equation}
By Remark~\ref{R2.2}, using the definition of the (Sobolev) generalized derivative of~$f$, it is found from~\eqref{eq4.38} that
\begin{equation}
\label{eq4.39}
f(x,t)-\varphi(x)=\int\limits_{0}^{t}D_{t}f(x,\tau)\,d\tau \qquad \hbox{for almost all } x \in \mathbb{R}^{n}.
\end{equation}
Next, by \eqref{eq4.39} and Remark~\ref{R2.2} we have, for any cube~$Q$,

\begin{gather}
\label{eq4.40}
\notag
\int\limits_{Q}|f(x,t)-\varphi(x)|\,dx \le \int\limits_{Q}\int\limits_{0}^{t}\biggl| D_{t}f(x,\tau)\,d\tau\biggr| \le \\
\notag
\le  C(C_{\gamma},Q)\|f|W^{1}_{1}(Q \times (0,t),\gamma)\| \to 0, \qquad t \to +0.
\end{gather}

The proof of the theorem is complete.

\begin{Def} \rm
Assume that a weight $\gamma \in A^{\rm loc}_{1}(\mathbb{R}^{n+1})$ and  $c_{1},c_{2} \geq 1$. By $Z=Z(\{\gamma_{k,m}\},c_{1},c_{2})$ we shall denote the linear space
of all functions $\varphi \in L^{\rm loc}_{1}(\mathbb{R}^{n})$ with finite norm (we set $E(\widetilde{Q}^{s}_{\alpha})\varphi:=E^{1}(\widetilde{Q}^{s}_{\alpha})\varphi$)
\begin{equation}
\label{eq4.41}
\|\varphi|Z\|:=\inf\limits_{T}\sum\limits_{s=1}^{\infty}\sum\limits_{\alpha \in \widetilde{A}^{s}}\widehat{\gamma}^{s}_{\alpha}E(\widetilde{Q}^{s}_{\alpha})\varphi + \sum\limits_{m \in \mathbb{Z}^{n}}\widehat{\gamma}_{0,m}\|\varphi|L_{1}(Q_{0,m})\|,
\end{equation}
where the infimum on the right of \eqref{eq4.41} is taken over  all tilings $T=\{T^{s}\}_{s=0}^{\infty}(c_{1}, c_{2})$
of the space~$\mathbb{R}^{n}$ that are \textit{admissible for the weight}~$\gamma$.
\end{Def}

The following main result of the present section is a direct corollary of Theorems \ref{Th4.1},~\ref{Th4.2} and the elementary estimate

$$
E(\widetilde{Q}^{s}_{\alpha}) \le \int\limits_{\widetilde{Q}^{s}_{\alpha}}|\varphi(x)-\varphi^{s}_{\alpha}|\,dx \le 2E(\widetilde{Q}^{s}_{\alpha}), \qquad s \in \mathbb{Z}_{+}, \alpha \in \widetilde{A}^{s}.
$$

\begin{Ca}
\label{C4}
Assume that a weight $\gamma \in A^{\rm loc}_{1}(\mathbb{R}^{n+1})$. Then there exist numbers
$c_{1} \geq q^{3}$, $c_{2} \geq q^{5}$ such that the operator $\operatorname{Tr}:W_{1}^{1}(\mathbb{R}^{n+1}_{+},\gamma) \to Z(\{\gamma_{k,m}\},c_{1},c_{2})$  is
continuous and there exists a~(nonlinear) continuous operator $\operatorname{Ext}: Z(\{\gamma_{k,m}\},c_{1},c_{2}) \to W_{1}^{1}(\mathbb{R}^{n+1}_{+},\gamma)$,
which is the right inverse of the operator $\operatorname{Tr}$.
\end{Ca}

\begin{Remark}\rm
From the proof of Theorems \ref{Th4.1}, \ref{Th4.2} it follows that for $c_{1} \geq q^{3}$ $c_{2} \geq q^{5}$ the space  $Z(\{\gamma_{k,m}\},c_{1},c_{2})$
is independent of the choice of constants $c_{1},c_{2}$, the corresponding norms being equivalent. Of course, the parameters $q^{3}$, $q^{5}$ may be fairly large.
But for us it is important that they are determined only from the sequence $\{\gamma_{k,m}\}$. Similarly, the space $Z(\{\gamma_{k,m}\},c_{1},c_{2})$
is independent of the choice of the parameter~$\lambda$ (which controls the expansion of the cubes $Q^{s}_{\alpha}$). Hence in what follows
the space $Z(\{\gamma_{k,m}\},c_{1},c_{2})$ will be denoted by $Z(\{\gamma_{k,m}\})$.
The following fairly subtle question is still open: find the constants $\sigma_{1}, \sigma_{2}$ such that for $c_{1} > \sigma_{1}$, $c_{2} > \sigma_{2}$, the corresponding
 norms in the space $Z(\{\gamma_{k,m}\})$ are equivalent, but for $c_{1} \le \sigma_{1}$ or $c_{2} \le \sigma_{2}$ the resulting norm is not equivalent
 to the norm of the space $Z(\{\gamma_{k,m}\})$. However, by author's opinion, this question plays no critical role for applications.
\end{Remark}

Let us establish some elementary properties of the space $Z(\{\gamma_{k,m}\})$.

\begin{Lm}
\label{Lm4.4}
Assume that a weight $\gamma \in A^{\rm loc}_{1}(\mathbb{R}^{n+1})$ and  $c_{1} \geq q^{3}$, $c_{2} \geq q^{5}$. Then, for the space $Z=Z(\{\gamma_{k,m}\})$,
we have the following \textbf{continuous} embeddings:
$$
\widetilde{B}^{1}(\mathbb{R}^{n},\{\gamma_{k,m}\}) \subset Z(\{\gamma_{k,m}\}) \subset L^{\rm loc}_{1}(\mathbb{R}^{n}).
$$
\end{Lm}

\textbf{The proof} of the continuity of the embedding $\widetilde{B}^{1}(\mathbb{R}^{n},\{\gamma_{k,m}\}) \subset Z(\{\gamma_{k,m}\})$ is clear. The second embedding  follows from
Corollary~\ref{C4},  Remark~\ref{R2.2} and the simple estimate

$$
\|\hbox{tr }|_{t=0} f|L_{1}(Q)\| \le  \|f|W^{1}_{1}(Q \times (0,1))\|,
$$
where $Q$ is a cube in the space $\mathbb{R}^{n}$.

\begin{Lm}
\label{Lm4.5}
Assume that a weight $\gamma \in A^{\rm loc}_{1}(\mathbb{R}^{n+1})$ and $c_{1} \geq q^{3}$, $c_{2} \geq q^{5}$. Then
the space $Z(\{\gamma_{k,m}\})=Z(\{\gamma_{k,m}\},c_{1},c_{2})$ is complete.
\end{Lm}

\textbf{The proof} follows from Corollary~\ref{C4} and the fact that the space $W_{1}^{1}(\mathbb{R}^{n+1}_{+},\gamma)$ is complete.

\begin{Remark}\rm
We claim that for $\gamma \equiv 1$ Gagliardo's result follows from Corollary~\ref{C4}.
The embedding $L_{1}(\mathbb{R}^{n}) \supset Z(\{\gamma_{k,m}\},c_{1},c_{2})$ with $c_{1} \geq q^{3}$, $c_{2} \geq q^{5}$ is clear.
To prove the converse embedding we note that
$$
\|\varphi|Z(\{\gamma_{k,m}\})\| \le \inf\limits_{\{l_{j}\}}\sum\limits_{j=0}^{\infty}\sum\limits_{m \in \mathbb{Z}^{n}}E(\widetilde{Q}_{l_{j},m})\varphi \le C \|\varphi|L_{1}(\mathbb{R}^{n})\|,
$$
where the infimum is taken over all sequences $\{l_{j}\}$ for which $l_{0}=0$ and the corresponding series is converging.
\end{Remark}


\begin{thebibliography}{99}

\bibitem{Tyu} A. I. Tyulenev Some new function spaces of variable smoothness, \textit{Sb. Mat.}, \textbf{224}, no. 6, 849-891 (2015).

\bibitem{Mi}
P. Mironescu and E. Russ. Traces of weighted sobolev spaces. Old and new. \textit{Nonlinear
Analysis: TMA}, 2015, 119, p.354-381.

\bibitem{Gal}
E. Gagliardo, Caratterizzazione delle trace sulla frontiera relative ad alcune classi di funzioni in $n$ variabili," \textit{Rend. Sem. Mat. Univ. Padova} 27, 284--305 (1957).

\bibitem{Tyu2}
A.~I.~Tyulenev, Description of traces of functions in the Sobolev space with a Muckenhoupt weight, \textit{Proc. Stekl. Inst. Math.} \textbf{284} , 280--295(2014).

\bibitem{Gin}
A.~S. ~Ginzburg Traces of functions from weighted classes, \textit{Izv. Vyssh. Uchebn. Zaved. Mat.}, 1984, no. 4, 61--64



\bibitem{Tyu3}
A.~I.~Tyulenev The Problem of Traces for Sobolev Spaces with Muckenhoupt-Type Weights, \textit{Mat. Zametki}, \textbf{94:5}, 720--732 (2013).


\bibitem{St}
E.~M.~Stein, \textit{Harmonic Analysis: Real-Variable methods, Orthogonality, and Oscillatory Integrals}, Princeton Univ. Press, Princeton, NJ, 1993.


\bibitem{Pe}
J. Peetre, A counterexample connected with Gagliardo's trace theorem, \textit{Comment. Math.} Special Issue 2 , 277---282(1979).


\bibitem{Sa}
I.~Mitsuo, Y. Sawano, Atomic decomposition for weighted Besov and Triebel--Lizorkin spaces, \textit{Math. Nachr.}, \textbf{285} , no.~1, 103--126(2012).

\bibitem{Ry}
V.S. Rychkov, Littlewood-Paley theory and function spaces with $A_{p}^{loc}$ - weights, \textit{Math. Nachr.}, \textbf{224}, 145--180(2001).



\bibitem{Bu}
V. I. Burenkov, \textit{Sobolev Spaces on Domains}, B. G. Teubner, Stuttgart, 1998.

\bibitem{Zi}
W.P.~Ziemer, \textit{Weakly differentiable functions: Sobolev Spaces and Functions of Bounded Variation}, Springer, New Orc, 1989.




\bibitem{Guz}
M.~de~Guzman, A covering lemma with application to differentiability of measures and singular integral operators, \textit{Studia Math.} \textbf{34}, 299--317 (1970).



\bibitem{Be}
O. V. Besov, To the Sobolev embedding theorem for the limiting exponent, \textit{Proc. Steklov Inst. Math.} \textbf{284}: 1, 81--96 (2014).

\bibitem{St2}
E. M. Stein, \textit{Singular integrals and differentiability properties of functions}, Princeton Univ. Press,
Princeton, New Jersey, 1970.

\end{thebibliography}
\end{document}